\documentclass[reqno, a4paper, 10pt]{amsart}

\usepackage[margin=1in]{geometry}
\usepackage[linesnumbered, noline, noend, ruled]{algorithm2e}
\usepackage[foot]{amsaddr}
\usepackage{amsfonts}
\usepackage{amsmath}
\usepackage{amssymb}
\usepackage{amsthm}
\usepackage[english]{babel}
\usepackage{bm}
\usepackage{booktabs}
\usepackage{braket}
\usepackage{float}
\usepackage[T1]{fontenc}
\usepackage{graphicx}
\usepackage[utf8]{inputenc}
\usepackage{nicefrac}
\usepackage[retain-explicit-plus, retain-zero-exponent, per-mode=symbol, output-exponent-marker=\text{e}, group-digits=false]{siunitx}
\usepackage[version=4]{mhchem}
\usepackage[table]{xcolor}
\usepackage{enumitem}
\usepackage{subcaption}

\usepackage{cite}
\usepackage{tikz}
\usepackage{hyperref}

\definecolor{myblue}{RGB}{0, 107, 164}
\definecolor{darkblue}{RGB}{0, 53, 82}
\definecolor{myorange}{RGB}{255, 128, 14}
\definecolor{darkorange}{RGB}{128, 63, 6}
\definecolor{grayone}{RGB}{171, 171, 171}
\definecolor{graytwo}{RGB}{89, 89, 89}
\definecolor{mybrown}{RGB}{200,82,0}

	\newtheorem*{Remark*}{Remark}

	\setenumerate{label=(\arabic*), ref=(\arabic*)}
	
	\renewcommand{\div}{\nabla \cdot}
	
	\newcommand{\opconstspace}{\thinspace}

	\newcommand{\R}{\mathbb{R}}

	\newcommand{\dmeas}[1]{\ \text{d}#1}
	
	\newcommand{\dual}[3]{\left\langle #1 , #2 \right\rangle_{#3}}
	\newcommand{\inner}[3]{\left( #1 , #2 \right)_{#3}}

	\makeatletter
	\DeclareOldFontCommand{\rm}{\normalfont\rmfamily}{\mathrm}
	\DeclareOldFontCommand{\sf}{\normalfont\sffamily}{\mathsf}
	\DeclareOldFontCommand{\tt}{\normalfont\ttfamily}{\mathtt}
	\DeclareOldFontCommand{\bf}{\normalfont\bfseries}{\mathbf}
	\DeclareOldFontCommand{\it}{\normalfont\itshape}{\mathit}
	\DeclareOldFontCommand{\sl}{\normalfont\slshape}{\@nomath\sl}
	\DeclareOldFontCommand{\sc}{\normalfont\scshape}{\@nomath\sc}
	\makeatother
	
	\numberwithin{equation}{section}
	
	\newcommand{\tablegray}{gray!25}
	
	\newcommand{\subin}{^\text{in}}
	\newcommand{\subwall}{^\text{wall}}
	\newcommand{\subout}{^\text{out}}
	\newcommand{\subreac}{^\text{reac}}
	\newcommand{\viscosity}{\mu}
	\newcommand{\density}{\rho}
	\newcommand{\hcapacity}{C_p}
	\newcommand{\conductivity}{\kappa}

	\newcommand{\normal}{\bm{n}}
	
	\newcommand{\velocity}{\bm{u}}
	\newcommand{\pressure}{p}
	\newcommand{\temperature}{T}
	\newcommand{\massfrac}{Y_k}
	\newcommand{\molefrac}{X_k}
	
	\newcommand{\conv}[1]{\chi^{#1}}
	
	\newcommand{\tempreacsource}{\dot{\omega}_{\temperature, \text{r}}}
	\newcommand{\tempdiffsource}{\dot{\omega}_{\temperature, \text{d}}}
	\newcommand{\specreacsource}{\dot{\omega}_{k}}
	\newcommand{\diffusioncoeff}{D_{k, \text{mix}}}
	\newcommand{\corrvelo}{\bm{V^c}}
	\newcommand{\pref}{p_{\text{ref}}}
	\newcommand{\ptot}{p_{\text{tot}}}
	\newcommand{\avmm}{M}
	\newcommand{\gasconst}{R}
	\newcommand{\specmm}{M_{k}}
	\newcommand{\spec}{\mathcal{M}_k}
	\newcommand{\stoich}{\nu}
	\newcommand{\kfj}{k_{\text{f}, j}}
	
	\newcommand{\keqj}{k_{\text{eq}, j}}
	\newcommand{\nspec}{N_{\text{s}}}
	\newcommand{\nreac}{N_{\text{r}}}
	\newcommand{\spechcap}{C_{p, k}}
	\newcommand{\specenth}{h_{k}}
	\newcommand{\specentr}{s_{k}}
	\newcommand{\specvis}{\viscosity_k}
	\newcommand{\speccond}{\conductivity_k}
	
	\newcommand{\eaj}{E_{\text{a}, j}}
	\newcommand{\patm}{\pressure_{\text{atm}}}
	\newcommand{\kf}{k_{\text{f}}}
	
	\newcommand{\keq}{k_{\text{eq}}}
	\newcommand{\nexp}{n}
	\newcommand{\logA}{\log(A)}
	\newcommand{\preexpj}{A_j}
	\newcommand{\preexp}{A}
	\newcommand{\ea}{E_{\text{a}}}
	\newcommand{\permeability}{K}
	\newcommand{\ihtc}{h_{\text{f,s}}}
	\newcommand{\scconst}{\chi^{\text{des}}}
	\newcommand{\enthalp}{h}
	\newcommand{\entrop}{s}
	\newcommand{\massfracvec}{\bm{Y^{\text{vec}}}}

	\newcommand{\grad}{\nabla}
	
	\newcommand{\integral}[1]{\int_{#1}}

	\newcommand{\costfunction}{J}

	\newcommand{\transposed}{^\top}
		
\usetikzlibrary{positioning, arrows, backgrounds, shadows, spy}

\begin{document}

\title{Optimal Control of the Sabatier Process in Microchannel Reactors}
\author{Sebastian Blauth$^{*,1,2}$}
\address{$^*$ Corresponding Author}
\address{$^1$ Fraunhofer ITWM, Kaiserslautern, Germany}
\email{\href{mailto:sebastian.blauth@itwm.fraunhofer.de}{sebastian.blauth@itwm.fraunhofer.de}}
\author{Christian Leith\"auser$^1$}
\email{\href{mailto:christian.leithaeuser@itwm.fraunhofer.de}{christian.leithaeuser@itwm.fraunhofer.de}}
\author{Ren\'e Pinnau$^2$}
\email{\href{mailto:pinnau@mathematik.uni-kl.de}{pinnau@mathematik.uni-kl.de}}
\address{$^2$ TU Kaiserslautern, Kaiserslautern, Germany}

\begin{center}
	This is a post-peer-review, pre-copyedit version of an article published in Journal of Engineering Mathematics. The final authenticated version is available online at:
	\url{https://dx.doi.org/10.1007/s10665-021-10134-2}
\end{center}

\begin{abstract}
	We consider the optimization of a chemical microchannel reactor by means of PDE constrained optimization techniques, using the example of the Sabatier reaction. To model the chemically reacting flow in the microchannels, we introduce a three- and a one-dimensional model. As these are given by strongly coupled and highly nonlinear systems of partial differential equations (PDEs), we present our software package cashocs which implements the adjoint approach and facilitates the numerical solution of the subsequent optimization problems. We solve a parameter identification problem numerically to determine necessary kinetic parameters for the models from experimental data given in the literature. The obtained results show excellent agreement to the measurements. Finally, we present two optimization problems for optimizing the reactor's product yield. First, we use a tracking-type cost functional to maximize the reactant conversion, keep the flow rate of the reactor fixed, and use its wall temperature as optimization variable. Second, we consider the wall temperature and the inlet gas velocity as optimization variables, use an objective functional for maximizing the flow rate in the reactor, and ensure the quality of the product by means of a state constraint. The results obtained from solving these problems numerically show great potential for improving the design of the microreactor.
	
	
\end{abstract}

\maketitle

\vspace{-0.5cm}
\section{Introduction}
\label{sec:introduction}

The Sabatier process, named after the French chemists Paul Sabatier and Jean-Baptiste Senderens who reported it in 1902 \cite{Senderens1902Nouvelles}, is given by the reversible exothermic reaction
\begin{equation}
	\label{eq:sabatier}
	\ce{CO2 + 4 H2 <=> CH4 + 2H2O}, \quad \Delta H^0 \approx \SI{-165}{\kilo\joule \per \mole}\quad (\text{at } \SI{25}{\celsius}).
\end{equation}
This reaction has been investigated, e.g., in the context of in-situ resource utilization on mars \cite{Carranza2005Microchannel, Hu2007Catalyst}, for life support systems on the ISS \cite{Samplatsky2012Development}, and it is also used for power-to-gas applications \cite{ElSibai2017Model, Falbo2018Kinetics}, and cogeneration systems \cite{Spazzafumo2018Cogeneration, Bailera2019Renewable}. An overview over these and various other applications can be found, e.g., in \cite{Vogt2019renaissance}. Microchannel geometries are particularly interesting for chemical reactions as the large specific surface area of the microchannels allows for high performance of catalytic reactions as well as precise temperature management by means of appropriate temperature control systems. Such reactors have already been investigated for the Sabatier reaction, e.g., in \cite{Brooks2007Methanation, Hu2007Catalyst, Engelbrecht2017Experimentation, Carranza2005Microchannel}.

The purpose of this paper is to investigate the optimization of such a microchannel reactor for the Sabatier reaction using techniques from PDE constrained optimization. Such methods are widely used in various physical applications, e.g., for the optimization of chemical reactions in conventional sized reactors \cite{vonSchwerin2000Process, Logist2011Robust, Benner2019Periodic}, the optimization of semiconductors \cite{Burger2003Fast, Hinze2002optimal}, glass cooling processes \cite{Thoemmes2002Numerical, Pinnau2004Optimal}, or for the optimal shape design of microchannel cooling systems \cite{Blauth2019Model, Blauth2020Shape}, aircrafts \cite{Schmidt2013Three, Schmidt2011Airfoil}, and polymer spin packs \cite{Leithaeuser2018Designing, Hohmann2019Shape}. The optimization of microchannel reactors has also been investigated previously using derivative free approaches only, e.g., in \cite{Tegrotenhuis2002Optimizing, Na2017Multi, Jeon2014Optimization, Jung2017Optimal}. To the best of our knowledge, the derivative based optimization of microchannel reactors with methods from PDE constrained optimization has not been considered in the literature so far.

Throughout this paper, we consider the same setting as in \cite{Engelbrecht2017Experimentation, Engelbrecht2017Carbon}, where a microchannel reactor for the Sabatier reaction is investigated by means of experiments and simulations. To model this reactor mathematically, we introduce the following two models. First, we present a three-dimensional model that contains all important physical and chemical effects of the reactor. Second, we derive a one-dimensional model from the first one using a homogenization procedure similar to \cite{Blauth2019Model}. As both models are given by strongly coupled and highly nonlinear systems of PDEs, we use our software package cashocs \cite{Blauth2020CASHOCS} for the numerical solution of the subsequent optimization problems. In particular, cashocs is used to determine the relevant kinetic parameters for both models by solving a parameter identification problem constrained by the one-dimensional reactor model. The obtained results show excellent agreement with the experimental results reported in \cite{Engelbrecht2017Experimentation, Engelbrecht2017Carbon}. Subsequently, a numerical comparison of both models shows that the one-dimensional model approximates the three-dimensional one very well, which validates our approach of using the one-dimensional model as the PDE constraint for the parameter identification.

Finally, we consider the following two optimization problems for the reactor. First, we use a cost functional based on the tracking of the \ce{CO2} conversion at the reactor outlet and treat the surrounding temperature of the reactor, i.e., its wall temperature, which can be influenced by means of appropriate temperature control systems, as optimization variable while keeping the inlet flow rate fixed. In the second case, we consider the inlet gas velocity of the reactor and its wall temperature as optimization variables and use a objective functional for maximizing the mass flow rate in the reactor. In this case, the quality of the product is ensured by use of a state constraint for the \ce{CO2} conversion. Both problems are solved numerically using our software cashocs and the obtained results show great potential for improving the design of microchannel reactors. 

This paper is structured as follows. We begin with the introduction of our mathematical models and an investigation of the Sabatier reaction in Section~\ref{sec:model_formulation}. Afterwards, in Section~\ref{sec:adpack}, we give a brief description of our software package cashocs which is used for the automated treatment of optimal control problems considered in the following sections. The necessary kinetic reaction parameters are determined in Section~\ref{sec:parameter_identification}, where we also compare both reactor models numerically. Finally, we investigate the numerical optimization of the reactor and discuss the potential for optimizing its design in Section~\ref{sec:optimal_control_reaction}.

\section{Model Formulation}
\label{sec:model_formulation}

We first give some preliminary notations regarding our setting. Afterwards we introduce a three- and a one-dimensional model for the reactor and investigate the behavior of the Sabatier reaction. We conclude this section by detailing the numerical solution procedure for both models.

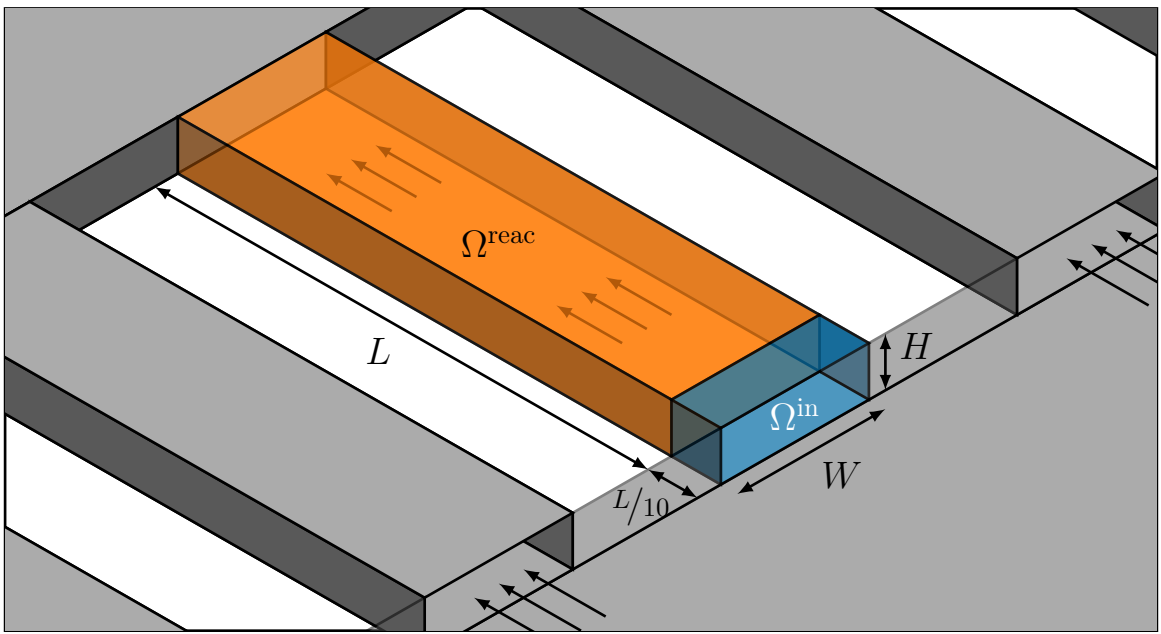
\begin{figure}[!b]
	\centering
%	\begin{tikzpicture}[scale=0.5, x={(1cm,0cm)}, y={(0cm, 1cm)}, z={(0.707107cm,0.707107cm)}, line width=1]
	\begin{tikzpicture}[scale=0.5, y={(0cm, 1cm)}, z={(-0.866cm,+0.5cm)}, x={(-0.866cm, -0.5cm)}, line width=1]
	
	\pgfmathsetmacro{\chanlen}{15}
	\pgfmathsetmacro{\chanwidth}{4.5}
	\pgfmathsetmacro{\chanheight}{1.5}
	\pgfmathsetmacro{\chandist}{4.5}
	\pgfmathsetmacro{\opac}{0.7}
	\pgfmathsetmacro{\transp}{0.5}
	
	%%% frame
	\draw (\chanwidth*2+\chandist*1.5+4.5,1,4.5) rectangle (-\chanwidth-\chandist*1.5-4.5,-1,10-4.5);
	\clip (\chanwidth*2+\chandist*1.5+4.5,1,4.5) rectangle (-\chanwidth-\chandist*1.5-4.5,-1,10-4.5);

	%%% collector
	%bottom
	\draw[fill=grayone] (-20,0,\chanlen) -- (15,0,\chanlen) -- (15,0,\chanlen*2) -- (-20,0,\chanlen*2) -- cycle;
	%top
	\draw[fill=grayone] (-20,\chanheight,\chanlen) -- (15,\chanheight,\chanlen) -- (15,\chanheight,\chanlen*2) -- (-20,\chanheight,\chanlen*2) -- cycle;
	%leftside
	\draw[fill=graytwo] (-20,0,\chanlen) -- (0,0,\chanlen) -- (0,\chanheight,\chanlen) -- (-20,\chanheight,\chanlen) -- cycle;
	%rightside
	\draw[fill=graytwo] (\chanwidth,0,\chanlen) -- (\chanwidth+10,0,\chanlen) -- (\chanwidth+10,\chanheight,\chanlen) -- (\chanwidth,\chanheight,\chanlen) -- cycle;

	%%% distributor
	%bottom
	\draw[fill=grayone] (-20,0,-\chanlen/10) -- (25,0,-\chanlen/10) -- (25,0,-\chanlen) -- (-20,0,-\chanlen) -- cycle;

	%%% main channel
	%back
%	\draw (0,0,\chanlen) rectangle (\chanwidth,\chanheight,\chanlen);
	%left
	\draw[fill=darkorange, opacity=\opac] (0,0,0) -- (0,0,\chanlen) -- (0,\chanheight,\chanlen) -- (0,\chanheight,0) -- cycle;
	%bottom
	\draw[fill=myorange, opacity=\opac] (0,0,0) -- (\chanwidth,0,0) -- (\chanwidth,0,\chanlen) -- (0,0,\chanlen) -- cycle;

	%%% Main Arrows Early
	\draw[-latex] (\chanwidth/2,\chanheight/2,\chanlen/5) -- (\chanwidth/2,\chanheight/2,\chanlen/3);
	\draw[-latex] (\chanwidth/3,\chanheight/2,\chanlen/5) -- (\chanwidth/3,\chanheight/2,\chanlen/3);
	\draw[-latex] (2*\chanwidth/3,\chanheight/2,\chanlen/5) -- (2*\chanwidth/3,\chanheight/2,\chanlen/3);
	
	%%% Main Arrows Late
	\draw[-latex] (\chanwidth/2,\chanheight/2,2/3*\chanlen) -- (\chanwidth/2,\chanheight/2,4/5*\chanlen);
	\draw[-latex] (\chanwidth/3,\chanheight/2,2/3*\chanlen) -- (\chanwidth/3,\chanheight/2,4/5*\chanlen);
	\draw[-latex] (2*\chanwidth/3,\chanheight/2,2/3*\chanlen) -- (2*\chanwidth/3,\chanheight/2,4/5*\chanlen);

	%front
%	\draw[fill=myorange, opacity=\opac] (0,0,0) rectangle (\chanwidth,\chanheight,0);
	%top
	\draw[fill=myorange, opacity=\opac] (0,\chanheight,0) -- (\chanwidth,\chanheight,0) -- (\chanwidth,\chanheight,\chanlen) -- (0,\chanheight,\chanlen) -- cycle;
	%right
	\draw[fill=darkorange, opacity=\opac] (\chanwidth,0,0) -- (\chanwidth,0,\chanlen) -- (\chanwidth,\chanheight,\chanlen) -- (\chanwidth,\chanheight,0) -- cycle;

	%%% artificial inlet
	%back
%	\draw[fill=myblue, opacity=\opac] (0,0,0) rectangle (\chanwidth, \chanheight,0);
	%left
	\draw[fill=darkblue, opacity=\opac] (0,0,-\chanlen/10) -- (0,0,0) -- (0,\chanheight,0) -- (0,\chanheight,-\chanlen/10) -- cycle;
	%bottom
	\draw[fill=myblue, opacity=\opac] (0,0,-\chanlen/10) -- (\chanwidth,0,-\chanlen/10) -- (\chanwidth,0,0) -- (0,0,0) -- cycle;
	%front
%	\draw[fill=myblue, opacity=\opac] (0,0,-\chanlen/10) rectangle (\chanwidth, \chanheight,-\chanlen/10);
	%right
	\draw[fill=darkblue, opacity=\opac] (\chanwidth,0,-\chanlen/10) -- (\chanwidth,0,0) -- (\chanwidth,\chanheight,0) -- (\chanwidth,\chanheight,-\chanlen/10) -- cycle;
	%top
	\draw[fill=myblue, opacity=\opac] (0,\chanheight,-\chanlen/10) -- (\chanwidth,\chanheight,-\chanlen/10) -- (\chanwidth,\chanheight,0) -- (0,\chanheight, 0) -- cycle;

	%%% right channel
	%left
	\draw[fill=graytwo] (\chanwidth+\chandist,0,-\chanlen/10) -- (\chanwidth+\chandist,0,\chanlen) -- (\chanwidth+\chandist,\chanheight,\chanlen) -- (\chanwidth+\chandist,\chanheight,-\chanlen/10) -- cycle;
	%bottom
	\draw[fill=grayone] (\chanwidth+\chandist,0,-\chanlen/10) -- (\chanwidth*2+\chandist,0,-\chanlen/10) -- (\chanwidth*2+\chandist,0,\chanlen) -- (\chanwidth+\chandist,0,\chanlen) -- cycle;
	%right
	\draw[fill=graytwo] (\chanwidth*2+\chandist,0,-\chanlen/10) -- (\chanwidth*2+\chandist,0,\chanlen) -- (\chanwidth*2+\chandist,\chanheight,\chanlen) -- (\chanwidth*2+\chandist,\chanheight,-\chanlen/10) -- cycle;
	%top
	\draw[fill=grayone] (\chandist+\chanwidth,\chanheight,-\chanlen/10) -- (\chandist+\chanwidth*2,\chanheight,-\chanlen/10) -- (\chandist+\chanwidth*2,\chanheight,\chanlen) -- (\chandist+\chanwidth,\chanheight,\chanlen) -- cycle;

	%%% right right channel
	%left
	\draw[fill=graytwo] (\chanwidth*2+\chandist*2,0,-\chanlen/10) -- (\chanwidth*2+\chandist*2,0,\chanlen) -- (\chanwidth*2+\chandist*2,\chanheight,\chanlen) -- (\chanwidth*2+\chandist*2,\chanheight,-\chanlen/10) -- cycle;
	%bottom
	\draw[fill=grayone] (\chanwidth*2+\chandist*2,0,-\chanlen/10) -- (\chanwidth*3+\chandist*2,0,-\chanlen/10) -- (\chanwidth*3+\chandist*2,0,\chanlen) -- (\chanwidth*2+\chandist*2,0,\chanlen) -- cycle;
	%right
	\draw[fill=graytwo] (\chanwidth*3+\chandist*2,0,-\chanlen/10) -- (\chanwidth*3+\chandist*2,0,\chanlen) -- (\chanwidth*3+\chandist*2,\chanheight,\chanlen) -- (\chanwidth*3+\chandist*2,\chanheight,-\chanlen/10) -- cycle;
	%top
	\draw[fill=grayone] (\chanwidth*2+\chandist*2,\chanheight,-\chanlen/10) -- (\chanwidth*3+\chandist*2,\chanheight,-\chanlen/10) -- (\chanwidth*3+\chandist*2,\chanheight,\chanlen) -- (\chanwidth*2+\chandist*2,\chanheight,\chanlen) -- cycle;

	%%% left channel
	%left
	\draw[fill=graytwo] (-\chanwidth-\chandist,0,-\chanlen/10) -- (-\chanwidth-\chandist,0,\chanlen) -- (-\chanwidth-\chandist,\chanheight,\chanlen) -- (-\chanwidth-\chandist,\chanheight,-\chanlen/10) -- cycle;
	%bottom
	\draw[fill=grayone] (-\chandist,0,-\chanlen/10) -- (-\chandist-\chanwidth,0,-\chanlen/10) -- (-\chandist-\chanwidth,0,\chanlen) -- (-\chandist,0,\chanlen) -- cycle;
	%top
	\draw[fill=grayone] (-\chandist,\chanheight,-\chanlen/10) -- (-\chandist-\chanwidth,\chanheight,-\chanlen/10) -- (-\chandist-\chanwidth,\chanheight,\chanlen) -- (-\chandist,\chanheight,\chanlen) -- cycle;
	%right
	\draw[fill=graytwo] (-\chandist,0,-\chanlen/10) -- (-\chandist,0,\chanlen) -- (-\chandist,\chanheight,\chanlen) -- (-\chandist,\chanheight,-\chanlen/10) -- cycle;

	%%% left left channel
	\draw[fill=grayone] (-\chanwidth*2-\chandist*2,0,-\chanlen/10) -- (-\chanwidth*2-\chandist*2,0,\chanlen) -- (-\chanwidth*2-\chandist*2,\chanheight,\chanlen) -- (-\chanwidth*2-\chandist*2,\chanheight,-\chanlen/10) -- cycle;
	\draw[fill=grayone] (-\chandist*2-\chanwidth,\chanheight,-\chanlen/10) -- (-\chandist*2-\chanwidth*2,\chanheight,-\chanlen/10) -- (-\chandist*2-\chanwidth*2,\chanheight,\chanlen) -- (-\chandist*2-\chanwidth,\chanheight,\chanlen) -- cycle;
	%right
	\draw[fill=graytwo] (-\chandist*2-\chanwidth,0,-\chanlen/10) -- (-\chandist*2-\chanwidth,0,\chanlen) -- (-\chandist*2-\chanwidth,\chanheight,\chanlen) -- (-\chandist*2-\chanwidth,\chanheight,-\chanlen/10) -- cycle;

	%%% Distributor
	%right right side
	\draw[fill=graytwo, opacity=\transp] (-\chandist*2-\chanwidth,0,-\chanlen/10) -- (-\chandist-\chanwidth,0,-\chanlen/10) -- (-\chandist-\chanwidth,\chanheight,-\chanlen/10) -- (-\chandist*2-\chanwidth,\chanheight,-\chanlen/10) -- cycle;
	% right side
	\draw[fill=graytwo, opacity=\transp] (-\chanwidth,0,-\chanlen/10) -- (0,0,-\chanlen/10) -- (0,\chanheight,-\chanlen/10) -- (-\chanwidth,\chanheight,-\chanlen/10) -- cycle;
	% left side
	\draw[fill=graytwo, opacity=\transp] (\chanwidth,0,-\chanlen/10) -- (\chanwidth+\chandist,0,-\chanlen/10) -- (\chanwidth+\chandist,\chanheight,-\chanlen/10) -- (\chanwidth,\chanheight,-\chanlen/10) -- cycle;
	% left left side
	\draw[fill=graytwo, opacity=\transp] (\chanwidth*2+\chandist,0,-\chanlen/10) -- (\chanwidth*2+\chandist*2,0,-\chanlen/10) -- (\chanwidth*2+\chandist*2,\chanheight,-\chanlen/10) -- (\chanwidth*2+\chandist,\chanheight,-\chanlen/10) -- cycle;
	%top
	%		\draw[fill=grayone, opacity=\transp] (-20,\chanheight,-\chanlen/10) -- (25,\chanheight,-\chanlen/10) -- (25,\chanheight,-\chanlen) -- (-20,\chanheight,-\chanlen) -- cycle;

	%%% Text descriptions
	\node at (\chanwidth/2, \chanheight/2,-\chanlen/10) {\LARGE \color{white} $\Omega\subin$};
	\node at (\chanwidth/2, \chanheight/2, \chanlen/2) {\LARGE $\Omega\subreac$};
	
	\draw [latex-latex] (0,-0,-\chanlen/10-0.5) -- node[below right] {\LARGE $W$} (\chanwidth,-0,-\chanlen/10-0.5);
	\draw [latex-latex] (\chanwidth+0.707,0,-\chanlen/10) -- node[below left=-0.125cm] {\LARGE $\nicefrac{L}{10}$} (\chanwidth+0.707,0,0);
	\draw [latex-latex] (\chanwidth+0.707,0,0) -- node[below left] {\LARGE $L$} (\chanwidth+0.707,0,\chanlen);
	\draw [latex-latex] (-0.5,0,-\chanlen/10) -- node[above right=-0.1cm and 0.05cm] {\LARGE $H$} (-0.5,\chanheight,-\chanlen/10);

	\draw[-latex] (-\chandist-\chanwidth/2,\chanheight/2,-2.5/15*\chanlen-\chanlen/10) -- (-\chandist-\chanwidth/2,\chanheight/2,-\chanlen/10);
	\draw[-latex] (-\chandist-\chanwidth/3,\chanheight/2,-2.5/15*\chanlen-\chanlen/10) -- (-\chandist-\chanwidth/3,\chanheight/2,-\chanlen/10);
	\draw[-latex] (-\chandist-2*\chanwidth/3,\chanheight/2,-2.5/15*\chanlen-\chanlen/10) -- (-\chandist-2*\chanwidth/3,\chanheight/2,-\chanlen/10);
	
	\draw[-latex] (\chanwidth+\chandist+\chanwidth/2,\chanheight/2,-2.5/15*\chanlen-\chanlen/10) -- (\chanwidth+\chandist+\chanwidth/2,\chanheight/2,-\chanlen/10);
	\draw[-latex] (\chanwidth+\chandist+\chanwidth/3,\chanheight/2,-2.5/15*\chanlen-\chanlen/10) -- (\chanwidth+\chandist+\chanwidth/3,\chanheight/2,-\chanlen/10);
	\draw[-latex] (\chanwidth+\chandist+2*\chanwidth/3,\chanheight/2,-2.5/15*\chanlen-\chanlen/10) -- (\chanwidth+\chandist+2*\chanwidth/3,\chanheight/2,-\chanlen/10);
	
	\draw[-latex] (\chanwidth*2+\chandist*2+\chanwidth/2,\chanheight/2,-2.5/15*\chanlen-\chanlen/10) -- (\chanwidth*2+\chandist*2+\chanwidth/2,\chanheight/2,-\chanlen/10);
	\draw[-latex] (\chanwidth*2+\chandist*2+\chanwidth/3,\chanheight/2,-2.5/15*\chanlen-\chanlen/10) -- (\chanwidth*2+\chandist*2+\chanwidth/3,\chanheight/2,-\chanlen/10);
	\draw[-latex] (\chanwidth*2+\chandist*2+2*\chanwidth/3,\chanheight/2,-2.5/15*\chanlen-\chanlen/10) -- (\chanwidth*2+\chandist*2+2*\chanwidth/3,\chanheight/2,-\chanlen/10);
	\end{tikzpicture}
	\caption{Schematic of the problem setting with inlet part $\Omega\subin$ (blue) and catalytic reaction part $\Omega\subreac$ (orange).}
	\label{fig:schematic}
\end{figure}

\subsection{Preliminary Notations}
\label{ssec:basic_notations}

For the modeling of the Sabatier process in a microchannel reactor we consider the experimental setting from \cite{Engelbrecht2017Experimentation, Engelbrecht2017Carbon}, where such a reactor is investigated by experiments and simulations. The reactor consists of 80 identical microchannels, each having a width of $W = \SI{450}{\micro\meter}$, a height of $H = \SI{150}{\micro\meter}$, and a length of $L = \SI{5}{\centi\meter}$ (cf. Table~\ref{tab:porous_parameters}). As in \cite{Engelbrecht2017Experimentation}, we assume that the flow is distributed uniformly between the channels so that we can model the behavior of the whole reactor by simulating a single channel, which is a well-established assumption in the literature \cite{Engelbrecht2017Experimentation, Chen2008Mathematical, Zeng2012Qualitative}. Additionally, we consider an inlet section, where we assume that no reaction occurs, with a length of $\nicefrac{L}{10}$ in front of the catalytic reaction part. Hence, throughout this paper the geometry of the entire reactor is denoted by $\Omega = (-\nicefrac{L}{10}, L) \times (0, W) \times (0, H)$ which is divided into the catalytic reactor domain $\Omega\subreac = (0, L) \times (0, W) \times (0, H)$ and the inlet domain $\Omega\subin = (-\nicefrac{L}{10}, 0) \times (0, W) \times (0, H)$, where we have no catalyst and, hence, can neglect the chemical reaction. The boundary $\Gamma = \partial\Omega$ of the channel is divided into three disjoint parts. The inlet $\Gamma\subin$ at $x=-\nicefrac{L}{10}$, where the fluid enters the domain, the wall boundary $\Gamma\subwall$, which bounds the domain and encloses the flow, and the outlet $\Gamma\subout$ at $x=L$, where the fluid leaves the domain. This setting is schematically depicted in Figure~\ref{fig:schematic}, where the geometry of the reactor is shown. For further details and schematics we refer to \cite{Engelbrecht2017Experimentation}. Finally, we remark that since the authors of \cite{Engelbrecht2017Experimentation, Engelbrecht2017Carbon} only considered stoichiometric inlet conditions, i.e., mole fractions of $\nicefrac{1}{5}$ for \ce{CO2} and $\nicefrac{4}{5}$ for \ce{H2}, we only consider these as inlet conditions throughout this paper.

\subsection{Mathematical Model}
\label{ssec:mathematical_model}

We assume that the chemically reacting gas mixture obeys the ideal gas law $\density = \frac{\ptot \avmm}{\gasconst \temperature}$, where $\ptot$ denotes the total pressure of the gas, $\avmm$ is its average molar mass, $\temperature$ its temperature, and $\gasconst$ is the universal gas constant. Further, we assume that the flow is weakly compressible, which implies that $\ptot$ can be written as $\ptot = \pref + \pressure$ with a constant reference pressure $\pref$ and pressure variations $\pressure$, so that $\pressure \ll \pref$ and $\grad \ptot = \grad \pressure$, as in \cite{Poinsot2005Theoretical}. Hence, we can approximate the ideal gas law by
\begin{equation}
	\label{eq:ideal_gas}
	\density = \frac{\pref \avmm}{\gasconst \temperature},
\end{equation}
as the contribution of the pressure variations towards the change in density of the gas is negligible. Moreover, for the derivation of a general model for chemically reacting fluid flow we assume that the reaction is described by the following system of $\nreac$ elementary and reversible chemical reactions between $\nspec$ chemical species, denoted by the symbol $\spec$ for $k=1,\dots, \nspec$,
\begin{equation}
	\label{eq:reaction_system}
	\sum_{k=1}^{\nspec} \stoich_{k, j}' \spec \ce{<=>} \sum_{k=1}^{\nspec} \stoich_{k, j}'' \spec, \quad j=1,\dots,\nreac,
\end{equation}
where $\stoich_{k, j}$ is the (non-negative) stoichiometric coefficient of species $k$ in reaction $j$. Note, that the reaction also affects the density of the gas mixture (cf.\ \eqref{eq:ideal_gas}) by changing its average molar mass $M$, which is defined as
\begin{equation*}
	\avmm = \left(\sum_{k=1}^{\nspec} \frac{\massfrac}{\specmm}\right)^{-1},
\end{equation*}
where $\massfrac$ is the mass fraction of species $k$, which are detailed below, and $\specmm$ is the average molar mass of species $k$.

We describe the chemically reacting flow in the channel by the following system of PDEs
\begin{equation}
	\label{eq:pde_system}
	\left\lbrace \quad
	\begin{alignedat}{2}
		\div \left(\density \velocity \right) &= 0 \quad && \text{ in } \Omega, \\
		\density \left( \velocity \cdot \grad \right) \velocity + \grad \pressure - \div \left( \viscosity \grad \velocity \right) - \grad \left( \nicefrac{\viscosity}{3}\ \div u \right) &= 0 \quad && \text{ in } \Omega, \\
		\density \hcapacity\ \velocity \cdot \grad \temperature - \div \left(\conductivity \grad \temperature \right) - \tempreacsource - \tempdiffsource &= 0 \quad && \text{ in } \Omega, \\
		\density\ \velocity \cdot \grad \massfrac + \div \left( \density \corrvelo \massfrac\right) - \div \left( \rho \diffusioncoeff \grad \massfrac \right) - \specreacsource &= 0 \quad &&\text{ in } \Omega, \quad k=1, \dots, \nspec,
	\end{alignedat}
	\right. 
\end{equation}
where $\velocity$ denotes the gas mixture's velocity, $\temperature$ its temperature, and $\grad \pressure$ corresponds to the gradient of the pressure (cf. the remarks above). Moreover, $\massfrac$ denotes the mass fraction of species $k$, which are combined to the vector $\massfracvec = [Y_1, \dots, Y_{\nspec}]\transposed$. Since the mass fractions sum to unity, we can get rid of, e.g., the last of the chemical species and compute its mass fraction via $Y_{\nspec} = 1 - \sum_{k=1}^{\nspec - 1} \massfrac$,
%\begin{equation*}
%Y_{\nspec} = 1 - \sum_{k=1}^{\nspec - 1} \massfrac,
%\end{equation*}
which decreases the computational complexity of the model. Further, $\density$, $\viscosity$, $\hcapacity$, and $\conductivity$ denote the gas mixture's density, viscosity, specific heat capacity, and thermal conductivity, respectively, and $\diffusioncoeff$ denotes the mixture-averaged diffusion coefficients. Note, that the equations of the PDE system \eqref{eq:pde_system} model the conservation of mass, momentum, enthalpy, and chemical species, respectively. 
%To model the diffusion of the chemical species we use the mixture-averaged diffusion coefficients $\diffusioncoeff$ and a correction velocity $\corrvelo$ that ensures the consistency between the mass and species conservation equations (cf. \cite{Poinsot2005Theoretical}). The term $\specreacsource$ models the conversion between the species due to the chemical reaction \eqref{eq:reaction_system}, and the terms $\tempreacsource$ and $\tempdiffsource$ model the change of temperature due to the reaction and molar diffusion of the species, respectively. 

The term $\specreacsource$ models the conversion between the species due to the chemical reaction \eqref{eq:reaction_system} and is given by
\begin{equation*}
	\specreacsource = \begin{cases}
	0 \quad &\text{ in } \Omega\subin, \\
	\specmm \sum_{j=1}^{\nreac} \stoich_{k, j} Q_j \quad &\text{ in } \Omega\subreac,
	\end{cases}
\end{equation*}
where we define $\stoich_{k, j} = \stoich_{k, j}'' - \stoich_{k, j}'$ and $\specmm$ is the molar mass of species $k$. Further, $Q_j$ is the rate of progress of reaction $j$, given by
\begin{equation}
	\label{eq:rate_of_progress}
	Q_j = \kfj \left( \prod_{k=1}^{\nspec} [\molefrac]^{\stoich_{k, j}'} - \frac{\prod_{k=1}^{\nspec} [\molefrac]^{\stoich_{k, j}''}}{\keqj}  \right),
\end{equation}
where $[\molefrac] = \frac{\density \massfrac}{\specmm}$ denotes the molar concentration of species $k$ and $\keqj$ is the equilibrium constant of reaction $j$, which is detailed in Appendix~\ref{app:constitutive_relations}. Further, $\kfj$ is the forward rate constant of reaction $j$, which is modeled by an Arrhenius law
\begin{equation*}
	\kfj = \preexpj \exp\left(- \frac{\eaj}{\gasconst \temperature} \right),
\end{equation*}
where $\preexpj$ is the so-called pre-exponential factor and $\eaj$ is the activation energy for reaction $j$. Note, that the reaction source term $\specreacsource$ vanishes in $\Omega\subin$ and is only active in $\Omega\subreac$, as discussed in Section~\ref{ssec:basic_notations}. Further, we remark that \eqref{eq:rate_of_progress} is only valid for elementary chemical reactions, hence our assumption for \eqref{eq:reaction_system}. 

The correction velocity $\corrvelo$, which is defined as
\begin{equation*}
	\corrvelo = \sum_{k=1}^{\nspec} \diffusioncoeff \grad \massfrac,
\end{equation*}
is used to ensure the consistency between the mass and species conservation equations for the mixture-averaged diffusion model, as is explained in detail, e.g., in \cite{Poinsot2005Theoretical, Kee2005Chemically}. The heat generated by the reaction is modeled through the source term $\tempreacsource$, which is given by
\begin{equation*}
	\tempreacsource = -\sum_{k=1}^{\nspec} \specenth(\temperature)\ \specreacsource,
\end{equation*}
where $\specenth$ is the pure species specific enthalpy (cf.\ Appendix~\ref{app:constitutive_relations}). Finally, the heat generated due to the molecular diffusion of the species is modeled by the term $\tempdiffsource$, which reads
\begin{equation*}
	\tempdiffsource = - \left( \density \sum_{k=1}^{\nspec} \spechcap \left( \massfrac \corrvelo - \diffusioncoeff \grad \massfrac \right) \right) \cdot \grad \temperature.
\end{equation*}

Note, that all transport parameters for the model depend on the fluid's chemical composition, represented by the mass fractions $\massfracvec$, and on the fluid's temperature $\temperature$ through constitutive relations. In particular, we have that
\begin{equation*}
	\density = \density(\massfracvec, \temperature), \quad \viscosity = \viscosity(\massfracvec, \temperature), \quad \hcapacity = \hcapacity(\massfracvec, \temperature), \quad \conductivity = \conductivity(\massfracvec, \temperature), \quad \diffusioncoeff = \diffusioncoeff(\massfracvec, \temperature),
\end{equation*}
which leads to a strong and highly nonlinear coupling between the individual equations of \eqref{eq:pde_system}. The constitutive relations detailing these dependencies are given in Appendix~\ref{app:constitutive_relations}.
%The constitutive equations detailing the transport parameters detailed in Appendix~\ref{app:constitutive_relations}, and can also be found in, e.g., \cite{Kee2005Chemically, Poinsot2005Theoretical}. Finally, we remark that we obtained the corresponding pure species parameters from \cite{McBride2002NASA, Svehla1995Transport}.

The system \eqref{eq:pde_system} is supplemented with the following boundary conditions. On the inlet $\Gamma\subin$ we use Dirichlet conditions to prescribe the ingoing velocity, temperature, and mass fractions of the gas mixture, i.e., we use
\begin{equation*}
%	\label{eq:bc_in}
	\left\lbrace \quad
	\begin{alignedat}{2}
		\velocity &= \velocity\subin \quad &&\text{ on } \Gamma\subin, \\
		\temperature &= \temperature\subin \quad &&\text{ on } \Gamma\subin, \\
		\massfrac &= \massfrac\subin \quad &&\text{ on } \Gamma\subin,\quad k=1,\dotsm,\nspec.
	\end{alignedat}
	\right.
\end{equation*}

On the wall boundary $\Gamma\subwall$ we use the usual no-slip condition for the velocity, which is valid as the Knudsen number of the flow is sufficiently small so that the continuum assumption holds. Additionally, we use a no-flux condition for the chemical species, which models that the species do not leave the reactor over the wall boundary, and a Dirichlet condition for the gas mixture's temperature. The latter models that the wall temperature of the reactor can be influenced through appropriate temperature control systems (cf. Section~\ref{ssec:realizability}). In summary, the boundary conditions on $\Gamma\subwall$ are given by
\begin{equation}
	\label{eq:bc_wall}
	\left\lbrace \quad
	\begin{alignedat}{2}
		\velocity &= 0 \quad &&\text{ on } \Gamma\subwall, \\
		\temperature &= \temperature\subwall \quad &&\text{ on } \Gamma\subwall, \\
		\density \left( \corrvelo \massfrac - \diffusioncoeff \grad \massfrac \right) \cdot \normal &= 0 \quad &&\text{ on } \Gamma\subwall,\quad k=1,\dots,\nspec,
	\end{alignedat}
	\right.
\end{equation}
where $\normal$ denotes the unit outer normal vector on $\Gamma$.

Finally, the unimpeded flow of the gas out of the reactor at the outlet $\Gamma\subout$ is modeled by a do-nothing condition for the momentum equation and homogeneous Neumann conditions for the temperature and chemical species, i.e.,
\begin{equation*}
%	\label{eq:bc_out}
	\left\lbrace \quad
	\begin{alignedat}{2}
		\viscosity \grad \velocity\ \normal + \frac{\viscosity}{3} \left(\div\velocity\right) \normal - \pressure\ \normal &= 0 \quad &&\text{ on } \Gamma\subout, \\
		\conductivity \grad \temperature \cdot \normal &= 0 \quad && \text{ on } \Gamma\subout, \\
		\density \left( \corrvelo \massfrac - \diffusioncoeff \grad \massfrac \right) \cdot \normal &= 0 \quad &&\text{ on } \Gamma\subout,\quad k=1,\dots,\nspec.
	\end{alignedat}
	\right.
\end{equation*}

\subsection{The Sabatier Reaction in Microchannel Reactors}
\label{ssec:sabatier_reaction}

The Sabatier reaction \eqref{eq:sabatier} is a reversible exothermic reaction used to convert \ce{CO2} and \ce{H2} into \ce{CH4} and \ce{H2O}. For most of its applications it is desirable that the reaction proceeds as far as possible to the right, so that ideally all of the \ce{CO2} would be consumed by the reaction. It is in this regard that we consider optimizing the microchannel reactor under investigation (cf. Section~\ref{sec:optimal_control_reaction}). As a measure for how far the reaction has already proceeded, we use the \ce{CO2} conversion, which is defined as
\begin{equation}
	\label{eq:def_conversion}
	\conv{\ce{CO2}} = 1 - \frac{Y_{\ce{CO2}}}{Y\subin_{\ce{CO2}}} = 1 - \frac{n_{\ce{CO2}}}{n\subin_{\ce{CO2}}},
\end{equation}
where $n\subin_{\ce{CO2}}$ denotes the molar amount of \ce{CO2} entering the domain. In this paper, we only evaluate $Y_{\ce{CO2}}$ and $n_{\ce{CO2}}$, which is the molar amount of \ce{CO2} in the gas mixture, at the outlet of the reactor (cf. Sections~\ref{sec:parameter_identification} and~\ref{sec:optimal_control_reaction}), so that $\conv{\ce{CO2}}$ describes the total \ce{CO2} conversion of the reactor.

Let us now investigate the behavior of the Sabatier reaction. First of all, we mention that there is another possible reaction between the reactants \ce{CO2} and \ce{H2}, namely the reverse water gas shift (RWGS) reaction, given by
\begin{equation*}
	\ce{CO2 + H2 <=> CO + H2O}, \quad \Delta H^0 \approx \SI{+41}{\kilo\joule \per \mole}\quad (\text{at } \SI{25}{\celsius}).
\end{equation*}
However, the RWGS reaction only becomes important for temperatures over approximately \SI{550}{\celsius} -- \SI{600}{\celsius} and shows negligible \ce{CO} production at lower temperatures \cite{Moioli2019Parametric, Mutschler2018CO2}. For these reasons, we proceed as in \cite{Moioli2019Parametric,Mutschler2018CO2} and restrict the temperatures under investigation to be below \SI{600}{\celsius} and neglect the RWGS reaction, which only introduces a limited error in our model.

The next point we need to address is that the Sabatier reaction is not an elementary reaction, which is required for the formula of the rate of progress \eqref{eq:rate_of_progress}. However, we decided not to break down the Sabatier reaction into a system of elementary reactions due to the following reasons. First, the Sabatier reaction is not fully understood yet from the viewpoint of elementary reaction mechanisms, with two possibilities proposed \cite{Baraj2016Reaction}. Second, using a system of elementary reactions would increase the computational complexity of the model significantly as this would introduce additional variables for the intermediate chemical species. For these reasons, we use the following modified formula for the rate of progress \eqref{eq:rate_of_progress} developed in \cite{Lunde1973Rates, Lunde1974Carbon, Lunde1974Modeling}
\begin{equation}
	\label{eq:empirical_rate_of_progress}
	Q = \kf \left( \left( \prod_{k=1}^{\nspec} [\molefrac]^{\stoich_{k, j}'} \right)^{\nexp} - \left( \frac{\prod_{k=1}^{\nspec} [\molefrac]^{\stoich_{k, j}''}}{\keq} \right)^{\nexp} \right),
\end{equation}
where $\nexp$ is an empirical exponent. Note, that \eqref{eq:empirical_rate_of_progress} models the Sabatier process using a single chemical reaction, so that we have $\nreac = 1$ and, hence, drop the index $j$ throughout the rest of the paper.
This approach has been used extensively in the literature to model the Sabatier reaction, e.g., in \cite{Engelbrecht2017Experimentation, Engelbrecht2017Carbon, Falbo2018Kinetics, Brooks2007Methanation, Moioli2019Parametric,Perez01Nov.2019Modeling, Moioli2019Model}. To model the forward rate constant $\kf$, we again use an Arrhenius law, i.e.,
\begin{equation}
	\label{eq:arrh}
	\kf = \preexp \exp\left( - \frac{\ea}{\gasconst \temperature} \right).
\end{equation}
It is straightforward to see that both formulations yield the same thermodynamic equilibrium as explained in \cite{Lunde1974Carbon}, which is a crucial feature of \eqref{eq:empirical_rate_of_progress}. Note, that even though this model for the Sabatier reaction is comparatively simple, it still covers the most important effects of the reactor. Hence, it is feasible that we consider this model to investigate the potential for improving the microchannel reactor. An investigation of more sophisticated reaction mechanisms or the inclusion of the RWGS reaction into the reactor model is beyond the scope of this paper and could be considered in future work.

\begin{figure}[!t]
	\centering
	\begin{subfigure}[t]{0.32\textwidth}
		\centering
		\includegraphics[width=\textwidth]{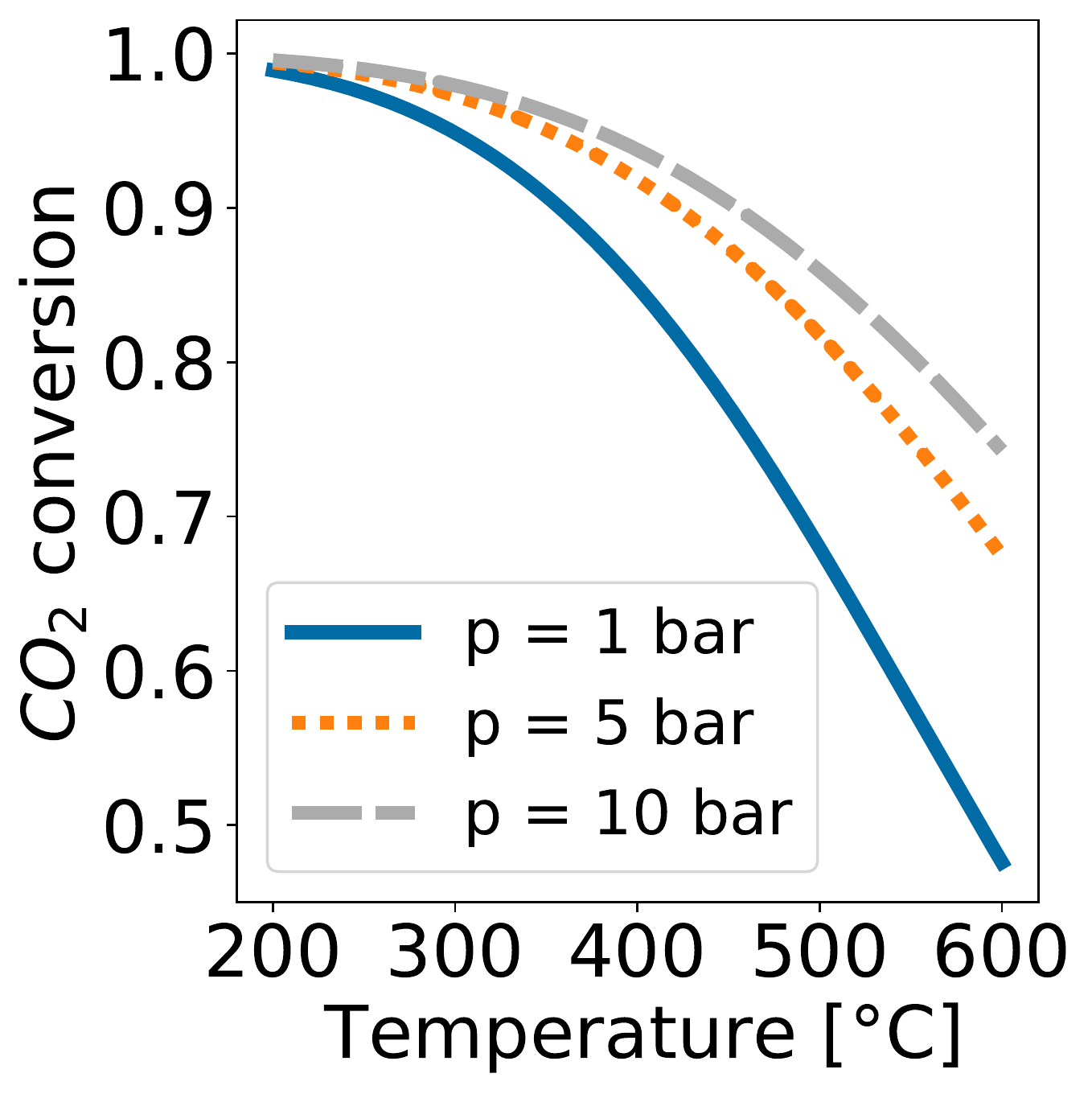}
		\caption{Equilibrium \ce{CO2} conversion for different pressures.}
		\label{sfig:sab_eq_conv}
	\end{subfigure}
	\hfil
	\begin{subfigure}[t]{0.32\textwidth}
		\centering
		\includegraphics[width=\textwidth]{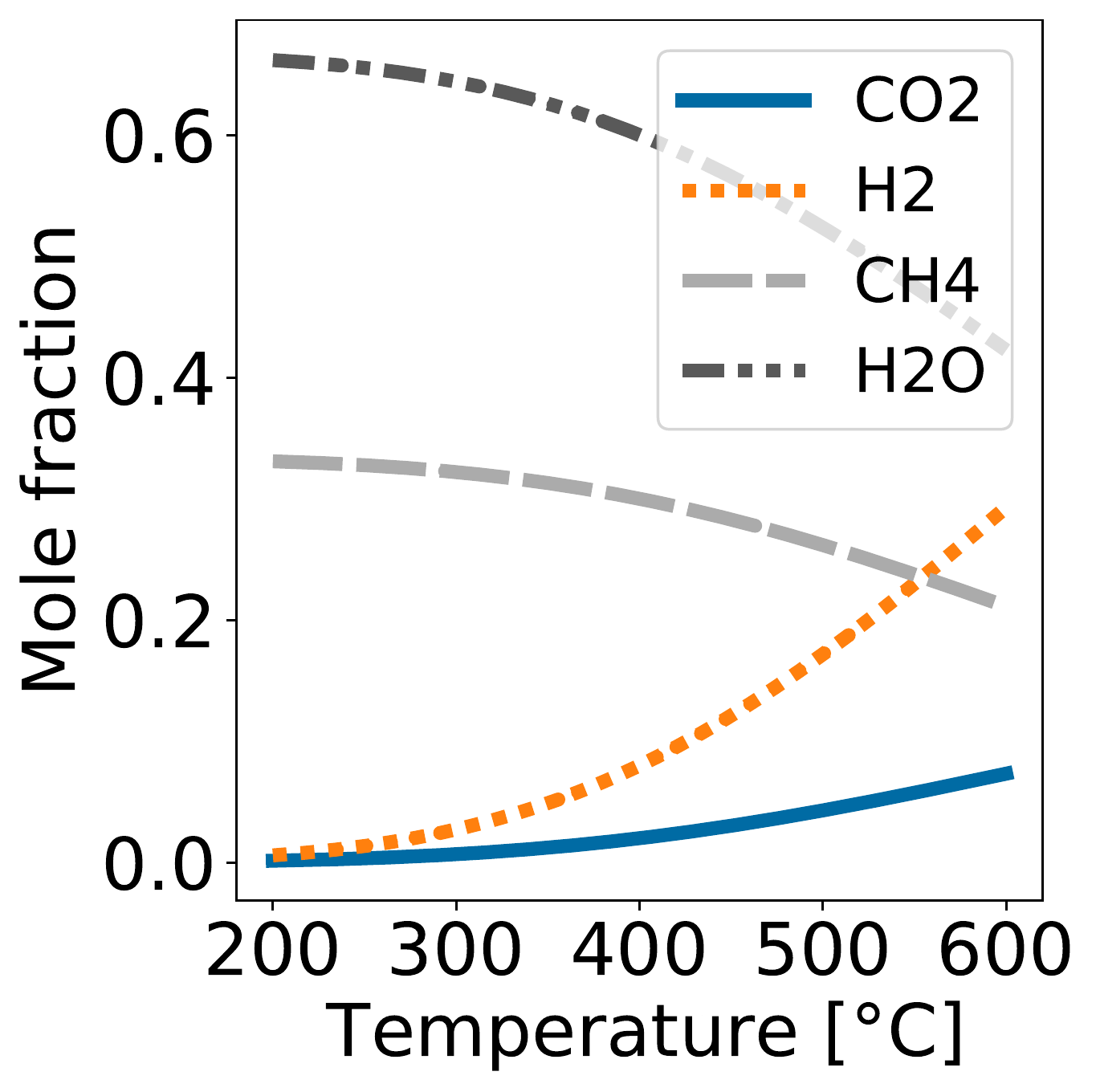}
		\caption{Equilibrium mole fractions at \SI{10}{\bar} pressure.}
		\label{sfig:sab_eq_mole}
	\end{subfigure}
	\hfil
	\begin{subfigure}[t]{0.32\textwidth}
		\centering
		\includegraphics[width=\textwidth]{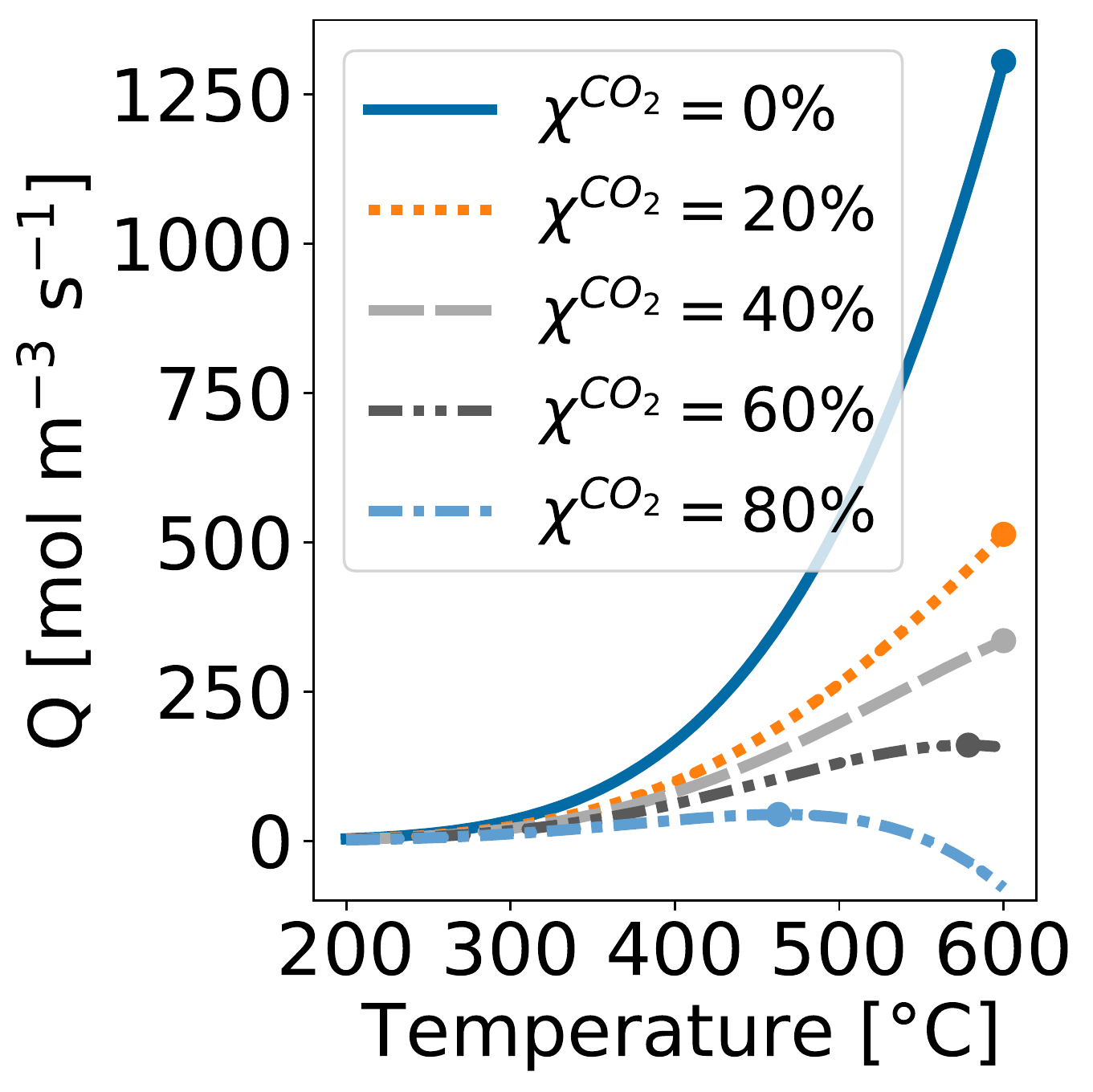}
		\caption{Rate of progress for different levels of \ce{CO2} conversion for \SI{10}{\bar} pressure.}
		\label{sfig:sab_reac_rates}
	\end{subfigure}
	\caption{Properties of the Sabatier reaction for stoichiometric inlet conditions.}
	\label{fig:sab_equilibrium}
\end{figure}

Let us now take a look at some important properties of the Sabatier reaction which are depicted in Figure~\ref{fig:sab_equilibrium} for stoichiometric inlet conditions (cf. Section~\ref{ssec:basic_notations}). In Figure~\ref{sfig:sab_eq_conv} the temperature dependence of the equilibrium \ce{CO2} conversion is shown for pressures of \SI{1}{\bar}, \SI{5}{\bar}, and \SI{10}{\bar}, which resemble the operating pressures investigated in \cite{Engelbrecht2017Experimentation, Engelbrecht2017Carbon}. We observe that the \ce{CO2} conversion decreases monotonically with temperature and increases monotonically with pressure. Both effects are direct consequences of Le Chatelier's principle \cite{Atkins2017Elements} applied to the Sabatier reaction. Hence, we only consider the case $\pref = \SI{10}{\bar}$ in this paper as the \ce{CO2} conversion is largest for this case. In Figure~\ref{sfig:sab_eq_mole} we can see the influence of the temperature on the equilibrium composition at the operating pressure of \SI{10}{\bar}. Again, we observe that the composition is more favorable for low temperatures of about \SI{200}{\celsius}, where it consists almost exclusively of the products \ce{CH4} and \ce{H2O}. The equilibrium composition becomes increasingly worse with higher temperatures as the amount of reactants increases considerably, mirroring our previous discussion regarding the equilibrium \ce{CO2} conversion. Finally, Figure~\ref{sfig:sab_reac_rates} shows the dependence of the rate of progress and, hence, of the reaction rate on the temperature at different levels of \ce{CO2} conversion, which we computed using the parameters determined later on (cf. Section~\ref{sec:parameter_identification}). We see that the rate of progress increases strongly with temperature for low levels of \ce{CO2} conversion. However, the greater the \ce{CO2} conversion, the lower the rate of progress becomes. For the highest levels of \ce{CO2} conversion, we observe that the maximum rate of progress, highlighted by the markers on the graphs, is obtained at progressively lower temperatures. Additionally, for sufficiently high levels of \ce{CO2} conversion and temperatures, the rate of progress even becomes negative and, thus, reverses the direction of the reaction. These are direct consequences of the thermodynamic equilibrium which pushes the reaction into the reverse direction under such conditions.

Regarding the optimization of the reactor, these characteristics lead to the following considerations. To obtain a high reaction rate, high temperatures are desirable, especially for low levels of \ce{CO2} conversion, as shown in Figure~\ref{sfig:sab_reac_rates}. However, if the temperature is kept high throughout the entire reactor, the thermodynamic limitations of the reaction severely constrain the maximum achievable \ce{CO2} conversion. On the other side, low temperatures are thermodynamically favorable, but the corresponding low reaction rate only leads to low levels of \ce{CO2} conversion. These considerations lead to the conclusion that, ideally, the reactor should have high temperatures near the inlet, where we have low levels of \ce{CO2} conversion, so that we have a high reaction rate initially. Subsequently, the temperature should decrease along the channel as this leads to more favorable equilibrium compositions and is also beneficial for the rate of progress at higher levels of \ce{CO2} conversion, which are confirmed in Section~\ref{sec:optimal_control_reaction}, where we investigate the optimization of the reactor. Finally, we note that the above deliberations regarding the influence of the temperature on the reaction is inherent in all exothermic reactions, and not specific to our choice of the Sabatier reaction, which is, again, a consequence of Le Chatelier's principle \cite{Atkins2017Elements}.

\subsection{A One-Dimensional Model of the Reactor}
\label{ssec:1D}

The PDE system \eqref{eq:pde_system} of Section~\ref{ssec:mathematical_model} completely models the behavior of chemically reacting fluid flow for a homogeneous reaction. However, a major drawback of the model is that it is three-dimensional, which makes the numerical simulation of the reactor comparatively costly and leads to a substantial bottleneck for the adjoint-based optimization of the reactor we consider in Sections~\ref{sec:parameter_identification} and~\ref{sec:optimal_control_reaction}. To remedy this, we now use ideas similar to \cite{Blauth2019Model} to derive a reduced one-dimensional model of the reactor which is significantly easier to solve numerically. 

As a first step, we interpret the geometry as a porous medium with porosity $1$. While the boundary $\Gamma\subwall$ is still present from this point of view due to the geometrical constraints, it now only acts as a mathematical boundary of the domain and does not interact physically with the fluid. To model the fluid velocity in such a porous medium, we use a Brinkman equation with slip boundary conditions, in analogy to \cite{Blauth2019Model} and the references therein, given by the system
\begin{equation*}
	\left\lbrace \quad
	\begin{alignedat}{2}
		\density \left( \velocity \cdot \grad \right) \velocity + \grad \pressure - \div \left(\viscosity \grad \velocity \right) - \grad \left( \nicefrac{\viscosity}{3}\ \div u \right) + \mu \permeability^{-1} \velocity &= 0 \quad && \text{ in } \Omega, \\
		\velocity \cdot \normal &= 0 \quad &&\text{ on } \Gamma\subwall, \\
		\mu\ \partial_{\normal} u \times \normal &= 0 \quad &&\text{ on } \Gamma\subwall,
	\end{alignedat}
	\right. 
\end{equation*}
where $\permeability$ denotes the permeability of the channel. Note, that the permeability is a property of the geometry only. It can be determined either analytically or numerically \cite{Bruus2007Theoretical}. Its value for our setting can be found in Table~\ref{tab:porous_parameters}. For the porous medium model, the above equations replace the momentum equation in \eqref{eq:pde_system} as well as the boundary condition for the velocity on $\Gamma\subwall$ in \eqref{eq:bc_wall}.

Additionally, we have to modify the temperature equation as we can no longer incorporate the effect of the wall temperature on the gas mixture via a Dirichlet boundary condition due to the previously mentioned reasons. Hence, we proceed analogously to \cite{Blauth2019Model} and use the following equation
\begin{equation*}
	\left\lbrace \quad
	\begin{alignedat}{2}
		\density \hcapacity\ \velocity \cdot \grad \temperature - \div \left(\conductivity \grad \temperature \right) + \ihtc \left( \temperature - \temperature\subwall \right) - \tempreacsource - \tempdiffsource &= 0 \quad && \text{ in } \Omega, \\
		\density \hcapacity\ u\cdot \normal - \conductivity \grad \temperature \cdot \normal &= 0 \quad && \text{ on } \Gamma\subwall, \\
	\end{alignedat}
	\right.
\end{equation*}
where $\ihtc$ is the so-called interfacial heat transfer coefficient which models the heat transfer between the gas mixture and the wall. We computed the value of $\ihtc$ numerically, as described in \cite{Blauth2019Model} and the references therein, and the result is shown in Table~\ref{tab:porous_parameters}. The no-flux boundary condition for the temperature again models that the temperature does not interact with the wall boundary. 

%These modifications yield a three-dimensional porous medium medium model that replaces the physical state variables that vary over the channel's cross section by averaged variables that are constant over the cross section. Due to the simple geometry of the channel and the fact that all boundary conditions on $\Gamma\subwall$ are chosen such that there is no physical interaction between the variables and the wall boundary, this implies that the model is quasi one-dimensional, i.e., all variables ($\velocity$, $\pressure$, $\temperature$, and $\massfrac$) only vary along the length of the channel. 

Averaging the system with the modifications described above over the channel cross section yields a one-dimensional model of the microchannel reactor, which is given by
%\begin{equation}
%	\label{eq:pde1d}
%	\left\lbrace \quad
%	\begin{alignedat}{2}
%		\div \left(\density \velocity \right) &= 0 \quad && \text{ in } \Omega, \\
%		%
%		\density \left( \velocity \cdot \grad \right) \velocity + \grad \pressure - \div \left(\viscosity \grad \velocity \right) - \grad \left( \nicefrac{\viscosity}{3}\ \div u \right) + \mu \permeability^{-1} \velocity &= 0 \quad && \text{ in } \Omega, \\
%		%
%		\density \hcapacity\ \velocity \cdot \grad \temperature - \div \left(\conductivity \grad \temperature \right) + \ihtc \left( \temperature - \temperature\subwall \right) - \tempreacsource - \tempdiffsource  &= 0 \quad && \text{ in } \Omega, \\
%		%
%		\density\ \velocity \cdot \grad \massfrac + \div \left( \density \corrvelo \massfrac\right) - \div \left( \rho \diffusioncoeff \grad \massfrac \right) - \specreacsource &= 0 \quad &&\text{ in } \Omega,\quad k=1, \dots, \nspec,
%	\end{alignedat}
%	\right. 
%\end{equation}
\begin{equation}
	\label{eq:pde1d}
	\left\lbrace \quad
	\begin{alignedat}{2}
		\partial_x \left(\density\opconstspace \velocity \right) &= 0 \quad && \text{ in } \Omega, \\
		\density \left( \velocity \opconstspace \partial_x \velocity \right) + \partial_x \pressure - \partial_x \left(\viscosity\opconstspace \partial_x \velocity \right) - \partial_x \left( \nicefrac{\viscosity}{3}\opconstspace \partial_x u \right) + \mu \permeability^{-1} \velocity &= 0 \quad && \text{ in } \Omega, \\
		\density \hcapacity\opconstspace \velocity\opconstspace \partial_x \temperature - \partial_x \left(\conductivity\opconstspace \partial_x \temperature \right) + \ihtc \left( \temperature - \temperature\subwall \right) - \tempreacsource - \tempdiffsource  &= 0 \quad && \text{ in } \Omega, \\
		\density\opconstspace \velocity\opconstspace \partial_x \massfrac + \partial_x \left( \density \corrvelo \massfrac\right) - \partial_x \left( \rho \diffusioncoeff\opconstspace \partial_x \massfrac \right) - \specreacsource &= 0 \quad &&\text{ in } \Omega,\quad k=1, \dots, \nspec,
	\end{alignedat}
	\right. 
\end{equation}
supplemented with the boundary conditions
%\begin{equation*}
%%	\label{eq:bc_in}
%	\left\lbrace \quad
%	\begin{alignedat}{2}
%		\velocity &= \velocity\subin \quad &&\text{ on } \Gamma\subin, \\
%		%
%		\temperature &= \temperature\subin \quad &&\text{ on } \Gamma\subin, \\
%		%
%		\massfrac &= \massfrac\subin \quad &&\text{ on } \Gamma\subin,\quad k=1,\dotsm,\nspec, \\
%		%
%		\viscosity \grad \velocity\ \normal + \frac{\viscosity}{3} \left(\div\velocity\right) \normal - \pressure\ \normal &= 0 \quad &&\text{ on } \Gamma\subout, \\
%		%
%		\conductivity \grad \temperature \cdot \normal &= 0 \quad && \text{ on } \Gamma\subout, \\
%		%
%		\density \left( \corrvelo \massfrac - \diffusioncoeff \grad \massfrac \right) \cdot \normal &= 0 \quad &&\text{ on } \Gamma\subout,\quad k=1,\dots,\nspec.
%	\end{alignedat}
%	\right.
%\end{equation*}
\begin{equation*}
	%	\label{eq:bc_in}
	\left\lbrace \quad
	\begin{alignedat}{2}
		\velocity &= \velocity\subin \quad &&\text{ on } \Gamma\subin, \\
		\temperature &= \temperature\subin \quad &&\text{ on } \Gamma\subin, \\
		\massfrac &= \massfrac\subin \quad &&\text{ on } \Gamma\subin,\quad k=1,\dotsm,\nspec, \\
		\viscosity\opconstspace \partial_x \velocity + \frac{\viscosity}{3}\opconstspace \partial_x \velocity - \pressure &= 0 \quad &&\text{ on } \Gamma\subout, \\
		\conductivity\opconstspace \partial_x \temperature &= 0 \quad && \text{ on } \Gamma\subout, \\
		\density \left( \corrvelo \massfrac - \diffusioncoeff\opconstspace \partial_x \massfrac \right) &= 0 \quad &&\text{ on } \Gamma\subout,\quad k=1,\dots,\nspec.
	\end{alignedat}
	\right.
\end{equation*}
Note, that the wall boundary is not present in \eqref{eq:pde1d} anymore, its influence is only modeled through the additional terms in the momentum and temperature equations, as discussed above. Additionally, for the sake of better readability we do not distinguish between the three-dimensional and the one-dimensional domain and denote both of them by $\Omega$ as it is obvious from the context to which we refer to. 

\begin{table}[!b]
	\centering
	{\footnotesize
		\rowcolors{2}{\tablegray}{white}
		\setlength{\tabcolsep}{1em}
		\begin{tabular}{r S}
			\toprule
			parameter [unit] & {value}\\
			\midrule
			length $L$ [\si{\meter}] & 5e-2 \\
			width $W$ [\si{\meter}] & 4.5e-4 \\
			height $H$ [\si{\meter}] & 1.5e-4 \\
			\midrule
			permeability $\permeability$ [\si{\cubic\meter}] & 1.48e-9 \\
			interfacial heat transfer coefficient $\ihtc$ [\si{\watt \per \kelvin \per \cubic \meter}] & 6.77e8 \\
			\bottomrule
		\end{tabular}
		\caption{Parameters for the reactor geometry and one-dimensional model \eqref{eq:pde1d}.}
		\label{tab:porous_parameters}
	}
\end{table}

Finally, we mention that for \eqref{eq:pde1d} only the form of the PDE system has changed slightly, all constitutive relations remain the same as described in Appendix~\ref{app:constitutive_relations} and Section~\ref{ssec:sabatier_reaction}. Furthermore, except for the kinetic parameters for the forward rate constant \eqref{eq:arrh}, all parameters for the model are obtained using the fits given in \cite{McBride2002NASA} and \cite{Svehla1995Transport} as discussed in Appendix~\ref{app:constitutive_relations}. We determine the missing parameters, namely the pre-exponential factor $\preexp$, the activation energy $\ea$, and the empirical exponent $\nexp$, by solving a parameter identification problem in Section~\ref{sec:parameter_identification}. Further, we compare the one-dimensional model \eqref{eq:pde1d} numerically to the three-dimensional one \eqref{eq:pde_system} in Section~\ref{ssec:comparison} after determining the aforementioned parameters. The corresponding results (cf. Figure~\ref{fig:error_analysis}) show that the one-dimensional model approximates the three-dimensional one very well, so that it is justified to use the former as model for the reactor throughout the rest of this paper.

\subsection{Numerical Solution of the Models}
\label{ssec:numerics}

To conclude this section, we briefly describe the methods used for solving the PDE systems \eqref{eq:pde_system} and \eqref{eq:pde1d}. First, for the one-dimensional model we discretize the corresponding interval with a uniform mesh consisting of 1001 nodes, corresponding to 1000 line segments. Second, for the three-dimensional model we also use a uniform mesh with 44040 nodes, corresponding to 178200 tetrahedrons. We use the finite element software FEniCS, version 2018.1 \cite{Alnes2015FEniCS, Logg2012Automated} to discretize both systems with a mixed finite element method using the following finite elements. The fluid's pressure, temperature, and mass fractions for \ce{CO2}, \ce{H2}, and \ce{CH4} are discretized with continuous, piecewise linear Lagrangian elements, and the fluid's velocity is discretized with continuous, piecewise quadratic Lagrangian elements. This yields the well-known Taylor-Hood finite element pair for velocity and pressure that is LBB stable for the saddle point structure of the continuity and momentum equations. As explained earlier, we calculate the mass fraction of \ce{H2O} via the relation $Y_{\ce{H2O}} = 1 - \left( Y_{\ce{CO2}} + Y_{\ce{H2}} + Y_{\ce{CH4}} \right)$. This discretization leads to a nonlinear system of equations which has to be solved. To do so, we use iterative methods which are described below and consider an initial guess given by zero velocity and pressure as well as a constant temperature and mass fractions determined by the respective inlet conditions for both models.

%Let us now briefly describe our discretization of the state and adjoint systems for all of the problems considered later on. We discretize the state system \eqref{eq:pde1d} using a mixed finite element method with continuous Lagrange elements. We decided to use the gas' velocity, pressure, temperature and the mass fractions of \ce{CO2}, \ce{H2}, and \ce{CH4} as variables for the finite element method. All other parameters, such as the density of the gas, are computed using the constitutive equations presented in Appendix~\ref{app:constitutive_relations} and Section~\ref{ssec:sabatier_reaction}. The mass fraction of \ce{H2O} is computed using the relation $Y_{\ce{H2O}} = 1 - \left( Y_{\ce{CO2}} + Y_{\ce{H2}} + Y_{\ce{CH4}} \right)$. We discretize all variables except the velocity with piecewise linear Lagrangian elements. The velocity, however, is discretized using piecewise quadratic Lagrangian elements. This is due to the fact that the continuity and momentum equations of \eqref{eq:pde1d} have a saddle-point like structure, such that an LBB stable pair of finite elements is needed for velocity and pressure. Note, that our choice corresponds to the Taylor-Hood finite element pair which is in fact LBB stable. Furthermore, for all optimization problems investigated we used a corresponding discretization of the adjoint states, i.e., the adjoint pressure, temperature, and mass fraction variables are discretized using piecewise linear Lagrangian elements and the adjoint velocity variable is discretized using piecewise quadratic Lagrangian elements. 

As the discretization described above leads to a system with 7006 degrees of freedom (DoF's) for the one-dimensional system \eqref{eq:pde1d}, we solve it monolithically using a damped Newton method with backtracking line search based on the natural monotonicity criterion described in \cite[Chapter 3.3]{Deuflhard2011Newton}. The arising linear systems for the Newton method are solved with the direct sparse linear solver MUMPS from the library PETSc \cite{Balay2020PETSc}.

The numerical solution of the resulting nonlinear system for the three-dimensional model \eqref{eq:pde_system} is more involved. As the system is considerably larger than the previously considered one, we do not use a monolithic approach for its solution. Instead we use a Picard-type fixed point iteration that consists of the following steps: first, we freeze the temperature and mass fractions, and solve the continuity and momentum equations to obtain values for the velocity and pressure. In the second step we freeze the velocity, pressure and temperature, and solve the species conservation equation for the mass fractions. The final step consists of freezing the velocity, pressure, and mass fractions, and solving the temperature equation. This is repeated until the residual of the original system reaches a relative tolerance of \num{1e-10}. Note, that in each of these steps a nonlinear system of equations has to be solved. This is done using the same damped Newton method as discussed above. In particular, the system of continuity and momentum equations has 176160 DoF's, the species conservation equation has 132120 DoF's, and the temperature equation has 44040 DoF's. As before, all resulting linear systems are solved using the solver MUMPS since they are sufficiently small to be solved by a direct method.

Since the size of the systems for the one-dimensional model is significantly smaller than the size of the ones for the three-dimensional model, the numerical solution of the former is also considerably faster. In particular, a solve of the one-dimensional model takes a few seconds, whereas it takes about an hour to solve the three-dimensional model, i.e., we get a speedup of over two orders of magnitude. This enables the fast solution of the optimization problems investigated subsequently in Sections~\ref{sec:parameter_identification} and~\ref{sec:optimal_control_reaction}.

%First, the continuity and momentum equations are solved together, giving an initial estimate of the flow field, while leaving $\massfrac$ and $\temperature$ fixed. Second, the species conservation equations are solved using the previously obtained velocity and a fixed temperature. Finally we solve the temperature equation, with velocity and mass fractions computed by the previous PDE solves. Then, we continue again with the continuity and momentum equations, until the residual reaches a relative tolerance of \num{1e-10}. Note, that each equation is solved is a nonlinear equation, as we only \qe{freeze} the nonlinearity induced by the coupling to the other equations, i.e., for each solve we only fixate the variables of the remaining equations and leave the variables belonging to the equation being solved free. To solve the nonlinear equations in the inner loop of the iteration, we again use a damped Newton method as described earlier. Finally, we remark that we also use the direct linear solver MUMPS for all arising linear systems, as the geometry is sufficiently simple even if it is three-dimensional.

\section{Automated Derivation of Adjoint Equations for the Numerical Solution of PDE Constrained Optimization Problems}
\label{sec:adpack}

The objective of this section is to briefly introduce our software package cashocs \cite{Blauth2020CASHOCS}, which implements the adjoint approach and is used to solve the PDE constrained optimization problems considered in Sections~\ref{sec:parameter_identification} and~\ref{sec:optimal_control_reaction}. To do so, we briefly recall the adjoint approach for PDE constrained optimization problems and then describe our implementation of this in cashocs, which mirrors this approach.

\subsection{Recapitulation of the Adjoint Approach}
\label{ssec:adjoint_approach}

Let us start by recalling the adjoint approach for PDE constrained optimization problems. For a detailed introduction to this topic we refer the reader, e.g., to the textbooks \cite{Hinze2009Optimization, Troeltzsch2010Optimal}. The general form of a PDE constrained optimization or optimal control problem is the following
\begin{equation}
	\label{eq:abstract_ocp}
	\min_{y, u}\ J(y, u) \quad \text{ s.t. } \quad e(y, u) = 0, \quad u\in U_\text{ad}, \quad \text{ and } \quad y\in Y_{\text{ad}},
\end{equation}
where $u \in U$ and $y \in Y$ are the so-called control and state variables that are part of appropriate Banach spaces $U$ and $Y$. Further, $J\colon Y \times U \to \R$ is the cost functional and $e\colon Y \times U \to Z^*$ is an operator between Banach spaces that models the PDE constraint, where $Z^*$ denotes the topological dual space of some Banach space $Z$. In particular, the so-called state equation $e(y, u) = 0$ is equivalent to
\begin{equation}
	\label{eq:abstract_state}
	\left\lbrace \quad
	\begin{aligned}
		&\text{Find } y \in Y \text{ such that }\\
		&\quad \dual{e(y, u)}{z}{Z^*, Z} = 0 \quad \text{ for all } z\in Z,
	\end{aligned}
	\right. 
\end{equation}
where $\dual{\varphi}{x}{B^*, B}$ denotes the duality pairing of $\varphi \in B^*$ and $x\in B$ for a Banach space $B$, which we also write as $\varphi[x]$. Note, that \eqref{eq:abstract_state} often corresponds to a variational form of a PDE. Finally, $U_\text{ad} \subset U$ and $Y_\text{ad} \subset Y$ are the non-empty and closed sets of admissible controls and states, respectively. These are used to model control and state constraints for the problem and can be treated by appropriate algorithms, such as projection methods for box constraints on the control variable or regularization techniques for state constraints \cite{Hinze2009Optimization, Troeltzsch2010Optimal}.

We assume that the state equation \eqref{eq:abstract_state} has a unique solution $y(u)$ for every $u\in U$ so that $e(y(u), u) = 0$ for all $u \in U$. This assumption allows us to define the so-called reduced cost functional $\hat{J} \colon U \to \R$ by
\begin{equation}
	\label{eq:abstract_reduced_cost_functional}
	\hat{J}(u) = J(y(u), u).
\end{equation}
In this way, we have formally eliminated the PDE constraint and can consider the following reduced optimization problem that is equivalent to \eqref{eq:abstract_ocp}
\begin{equation*}
	\min_{u}\ \hat{J}(u) \quad \text{ s.t. } \quad u\in U_\text{ad} \quad \text{ and } \quad y(u) \in Y_{\text{ad}}.
\end{equation*}
Our goal is to efficiently compute the gradient of the reduced cost functional $\hat{J}'(u)$ which is to be used as part of derivative based optimization algorithms for solving \eqref{eq:abstract_ocp}. To do so, we assume that $J$ and $e$ are continuously Fr\'echet differentiable, and that $e_y(y(u), u)$, i.e., the Fr\'echet derivative of $e$ w.r.t.\ $y$, is continuously invertible, so that the implicit function theorem \cite{Hinze2009Optimization} ensures that $y(u)$ is continuously differentiable in a neighborhood of a solution of \eqref{eq:abstract_state}.
%Under these assumptions, $y'(u)$ can be calculated using
%\begin{equation*}
%	0 = \frac{d}{d u} e(y(u), u) = e_y(y(u), u) [y'(u)] + e_u(y(u), u) \quad \Rightarrow  \quad y'(u) = e_y^{-1}(y(u), u) [e_u(y(u), u)].
%\end{equation*}
%Now, using the chain rule allows us to compute the sensitivities of the reduced cost functional as
%\begin{equation*}
%	\dual{\hat{J}'(u)}{h}{U^*, U} = \dual{J_y(y(u), u)}{y'(u)[h]}{Y^*, Y} + \dual{J_u(y(u), u)}{h}{U^*, U}.
%\end{equation*}
%However, as is explained in \cite{Hinze2009Optimization}, it is very costly to compute $y'(u)$. For this reason, we proceed differently and introduce an adjoint variable $p \in Z$ for the PDE constraint and set up a Lagrangian function as follows
To calculate the gradient of the reduced cost functional, we now derive the necessary adjoint and gradient equations using a Lagrange approach. We introduce an adjoint variable $p \in Z$ for the state equation \eqref{eq:abstract_state} and set up a Lagrangian function $L\colon Y\times U \times Z \to \R$ corresponding to \eqref{eq:abstract_ocp} as follows
\begin{equation}
	\label{eq:def_lagrangian}
	L(y, u, p) = J(y, u) + \dual{e(y, u)}{p}{Z^*, Z}.
\end{equation}
If we insert $y = y(u)$ into the Lagrangian, the PDE constraint vanishes and we have
\begin{equation*}
	L(y(u), u, p) = J(y(u), u) + \dual{e(y(u), u)}{p}{Z^*, Z} = \hat{J}(u) \quad \text{ for all } p \in Z.
\end{equation*}
Hence, differentiating the Lagrangian at $(y(u), u, p)$ w.r.t.\ $u$ yields
\begin{equation}
	\label{eq:gradient_computation_help}
	\dual{\hat{J}'(u)}{h}{U^*, U} = \dual{\frac{d}{du} L(y(u), u, p)}{h}{U^*, U} = \dual{L_y(y(u), u, p)}{y'(u)[h]}{Y^*, Y} + \dual{L_u(y(u), u, p)}{h}{U^*, U}.
\end{equation}
The main idea of the adjoint approach is now to choose $p = p(u) \in Z$ so that $L_y(y(u), u, p) = 0$, i.e.,
\begin{equation*}
	\dual{L_y(y(u), u, p)}{q}{Y^*, Y} = 0 \quad \text{ for all } q\in Y.
\end{equation*}
This equation can be rewritten as
\begin{align*}
	0 = \dual{L_y(y(u), u, p)}{z}{Y^*, Y} &= \dual{J_y(y(u), u)}{z}{Y^*, Y} + \dual{e_y(y(u), u)[z]}{p}{Z^*, Z} \\
	&= \dual{J_y(y(u), u) + e^*_y(y(u), u)[p]}{z}{Y^*, Y},
\end{align*}
and we note that the use of the adjoint operator $e^*_y$ gives the adjoint approach its name. In particular, the above equation can be interpreted as
\begin{equation}
	\label{eq:abstract_adjoint}
	\left\lbrace \quad
	\begin{aligned}
		&\text{Find } p \in Z \text{ such that } \\
		& \quad \dual{e_y^*(y(u), u)[p]}{z}{Y^*, Y} = - \dual{J_y(y(u), u)}{z}{Y^*, Y}.
	\end{aligned}
	\right.
\end{equation}
As for the state equation, we assume that the adjoint equation \eqref{eq:abstract_adjoint} is well-posed, and the corresponding solution $p = p(u)$ is called the adjoint state w.r.t.\ $u$. Furthermore, it is worth noticing that the adjoint equation \eqref{eq:abstract_adjoint} is always a linear equation, usually simplifying its numerical solution significantly compared to the possibly nonlinear state equation. 
Inserting $p=p(u)$ into \eqref{eq:gradient_computation_help} reveals
\begin{equation}
	\label{eq:abstract_riesz}
	\begin{aligned}
		\dual{\hat{J}'(u)}{h}{U^*, U} &= \dual{L_y(y(u), u, p(u))}{y'(u)[h]}{Y^*, Y} + \dual{L_u(y(u), u, p(u))}{h}{U^*, U} \\
		&= \dual{L_u(y(u), u, p(u))}{h}{U^*, U} \\
		&= \dual{J_u(y(u), u)}{h}{U^*, U} + \dual{e_u(y(u), u) [h]}{p(u)}{Z^*, Z} \\
		&= \dual{J_u(y(u), u) + e_u^*(y(u), u) [p(u)]}{h}{U^*, U}.
	\end{aligned}
\end{equation}
However, \eqref{eq:abstract_riesz} only specifies the derivative of the reduced cost functional, i.e., an element of $U^*$. In case that $U$ is a Hilbert space, we can use the Riesz representation theorem to identify this with the gradient of the reduced cost functional, which is then given by
\begin{equation*}
	\hat{J}'(u) = J_u(y(u), u) + e_u^*(y(u), u) [p(u)] \in U.
\end{equation*}

In summary, to compute the gradient of the reduced cost functional \eqref{eq:abstract_reduced_cost_functional} for a given $u\in U$, we solve the state equation \eqref{eq:abstract_state} to obtain the state variable $y(u)$ and then solve the adjoint equation \eqref{eq:abstract_adjoint} to determine the adjoint variable $p(u)$. If $U$ is a Hilbert space, we compute the gradient $\hat{J}'(u)$ using a Riesz projection as discussed above.

\subsection{Automated Derivation of Adjoint and Gradient Equations}
\label{ssec:ad}

From our recapitulation of the adjoint approach we see that to derive the adjoint and gradient equations \eqref{eq:abstract_adjoint} and \eqref{eq:abstract_riesz}, we need to differentiate the cost functional and state equation w.r.t.\ the state $y$ and the control $u$ and to calculate the corresponding adjoint operators. For many optimization problems with complex, coupled, and nonlinear PDE constraints, it is often impossible to verify all of the necessary assumptions of the adjoint approach. Moreover, even assuming that all objects are sufficiently smooth and calculating the adjoint and gradient equations in a formal manner is not feasible for very complex PDE constraints as it involves extremely tedious and error-prone calculations \cite{Farrell2013Automated}. Our reactor models from Section~\ref{sec:model_formulation} certainly belong to this category of PDE constraint. For these reasons, we use our software package cashocs \cite{Blauth2020CASHOCS}, which offers a numerical implementation of the continuous adjoint approach based on the finite element software FEniCS \cite{Alnes2015FEniCS, Logg2012Automated}. In particular, it uses the built-in automatic differentiation capabilities of FEniCS to derive the corresponding adjoint and gradient equations. 

In FEniCS, the cost functional and the PDE constraint of many optimization problems can be represented as variational forms in the so-called Unified Form Language (UFL). These forms can be differentiated automatically and symbolically with respect to their arguments, as explained in \cite[Chapter 17]{Logg2012Automated}.
%We use this feature to automatically derive the adjoint and gradient equations as discussed in the previous section. Hence, in order to derive the adjoint and gradient equations, the user has to specify the forms corresponding to the state system and the cost functional, as well as to declare, which variables are the state and which are the control variables. 
To derive the corresponding variational forms of the adjoint and gradient equations, cashocs sets up a UFL form for the Lagrangian as in \eqref{eq:def_lagrangian} which is subsequently differentiated symbolically w.r.t.\ the state and control variables \cite{Blauth2020CASHOCS}. Finally, the gradient $g \in U$ of the reduced cost functional w.r.t.\ some Hilbert space $U$ is computed numerically by solving the Riesz problem 
\begin{equation*}
	\left\lbrace\quad
	\begin{aligned}
		&\text{Find } g \in U \text{ such that } \\
		&\quad \inner{g}{h}{U} = \inner{\hat{J}'(u)}{h}{U} = \inner{L_u(y(u), u, p(u))}{h}{U} \quad \text{ for all } h \in U,
	\end{aligned}
	\right.
\end{equation*}
where $\inner{\cdot}{\cdot}{U}$ denotes the scalar product of $U$. Throughout this paper we only consider the case $U=L^2(\Omega)$ or $U=L^2(\Gamma\subout)$, in particular, we use the $L^2$ scalar product to identify the gradients of the cost functionals. 

We remark that the two differentiation operations described above are the only places where automatic differentiation is used in cashocs, the rest of the package, which manages the optimization problem, as well as the optimization algorithms do not rely on automatic differentiation or pre-existing implementations. Note, that if we discretize the state problem with a Galerkin method that uses the continuous variational formulation, as we do for the models considered in this paper, our software cashocs computes the corresponding continuous variational forms of the adjoint and gradient equations, which are only discretized at a later stage. In this case, cashocs consistently discretizes the continuous adjoint approach we recalled previously \cite{Blauth2020CASHOCS}. Finally, we note that there are software packages that use similar ideas to derive adjoint equations utilizing the AD capabilities of the UFL, e.g., dolfin-adjoint \cite{Mitusch2019dolfin} or firedrake \cite{Rathgeber2017Firedrake, Ham2019Automated}. However, we do not use these packages as our own approach can be tailored better to our needs and gives us complete control over the optimization algorithms used for solving the corresponding problems.

\section{Identification of Kinetic Reaction Parameters}
\label{sec:parameter_identification}

In this section, we determine the kinetic parameters needed to complete our models of the reactor (cf. Section~\ref{ssec:1D}) from the experimental results reported in \cite{Engelbrecht2017Experimentation, Engelbrecht2017Carbon}. We introduce a parameter identification problem based on our one-dimensional reactor model which is solved numerically utilizing our software package cashocs described in Section~\ref{sec:adpack}. Finally, we compare our two reactor models from Section~\ref{sec:model_formulation} numerically and see that both yield nearly identical results.

\subsection{Description of the Parameter Identification Problem}

In \cite{Engelbrecht2017Experimentation, Engelbrecht2017Carbon}, the authors consider a different model for the reactor with the following major variations to our models from Section~\ref{sec:model_formulation}. They consider a heterogeneous reaction only occurring in a porous catalyst located close to the channel wall, whereas we consider a homogeneous reaction. Moreover, they assume that the reaction starts immediately at the inlet of the computational domain, whereas we include the inlet section $\Omega\subin$where no reaction occurs. For these reasons, we cannot use the kinetic parameters reported in \cite{Engelbrecht2017Experimentation, Engelbrecht2017Carbon}, but have to determine appropriate ones for our models ourselves.
%values for the forward rate constant $\kf$ in \eqref{eq:arrh}, namely the activation energy $\ea$, the pre-exponential factor $\preexp$, and the empirical exponent $\nexp$, for our models ourselves.

As remarked in Section~\ref{ssec:sabatier_reaction}, we restrict our investigation to the case of \SI{10}{\bar} for the operating pressure of the reactor. For this setting, a total of 21 experiments, considering 7 different reactor temperatures $\temperature\subwall$ (\SI{250}{\celsius} -- \SI{400}{\celsius} with increments of \SI{25}{\celsius}) and 3 different inlet flow rates (\SI{50}{\milli\liter \per \minute}, \SI{100}{\milli\liter \per \minute}, and \SI{150}{\milli\liter \per \minute}), were carried out in \cite{Engelbrecht2017Experimentation, Engelbrecht2017Carbon}. Note, that here and throughout the rest of this paper, when we specify flow rates of the gas, we always assume normal conditions given by a pressure of \SI{1}{atm} and a temperature of \SI{273.15}{\kelvin} so that we can compare them regardless of the physical conditions. We remark that the experimental results of \cite{Engelbrecht2017Experimentation, Engelbrecht2017Carbon} are given by measurements of the achieved \ce{CO2} conversion, and that we extracted the corresponding numerical values using the software Webplot Digitizer \cite{Rohatgi2019Webplotdigitizer}.

To determine the kinetic parameters, we consider the following parameter identification problem
\begin{equation}
	\label{eq:parameter_identification_problem}
	\min_{y, u_c}\ J(y, u_c) = \sum_{l=1}^{21} \frac{1}{2} \integral{\Gamma\subout} \left( \conv{\ce{CO2}}_{\text{sim}, l} - \conv{\ce{CO2}}_{\text{exp}, l} \right)^2 \dmeas{s} \quad \text{ s.t. } \quad e(y, u_c) = 0 \quad \text{ and } \quad u_c \in U_\text{ad},
\end{equation}
where the state variables are combined in the vector $y = [\pressure, \velocity, \temperature, \massfracvec]\transposed$ and the control variables are combined in the vector $u_c = [\ea, \logA, \nexp]\transposed$. Note, that we use the logarithm of the pre-exponential factor, i.e., $\logA$, instead of $\preexp$ as control variable since this ensures a better scaling of the parameters and, additionally, guarantees that the computed pre-exponential factor is positive. Moreover, we assume that $u_c$ is constant in $\Omega$, i.e., that the kinetic parameters do not vary spatially. This yields a finite dimensional optimization problem, so that we do not require additional regularization for the cost functional. Furthermore, $\conv{\ce{CO2}}_{\text{exp}, l}$ denotes the \ce{CO2} conversion measured in experiment $l$, and $\conv{\ce{CO2}}_{\text{sim}, l}$ is the \ce{CO2} conversion of the $l$-th simulation, which are calculated via \eqref{eq:def_conversion}. The operator for the state system $e(y, u_c)$ is given by
\begin{equation*}
	e(y, u_c) = [e_1(y, u_c), \dots, e_{21}(y, u_c)]\transposed,
\end{equation*}
where $e_l(y, u_c) = 0$ is the weak form of the state system \eqref{eq:pde1d} with appropriate values for $\velocity\subin$ and $\temperature\subwall$ corresponding to the $l$-th experiment, as discussed above. Finally, the set of admissible controls is given by
\begin{equation*}
	U_\text{ad} = \Set{[\ea, \logA, \nexp]\transposed \in \R^3 | \nexp \geq 0},
\end{equation*}
which models that the empirical exponent is supposed to be non-negative. Note, that the cost functional in \eqref{eq:parameter_identification_problem} corresponds to a least-squares problem for fitting the \ce{CO2} conversion to the experimental results. 
%Finally, we remark that we did not use any regularization term for this cost functional as the control variables are scalar quantities and we were, nevertheless, able to solve \eqref{eq:parameter_identification_problem} very well numerically.

\subsection{Numerical Results}
\label{ssec:results_identification}

\begin{figure}[!b]
	\centering
	\begin{subfigure}[t]{0.49\textwidth}
		\centering
		\includegraphics[width=0.9\textwidth]{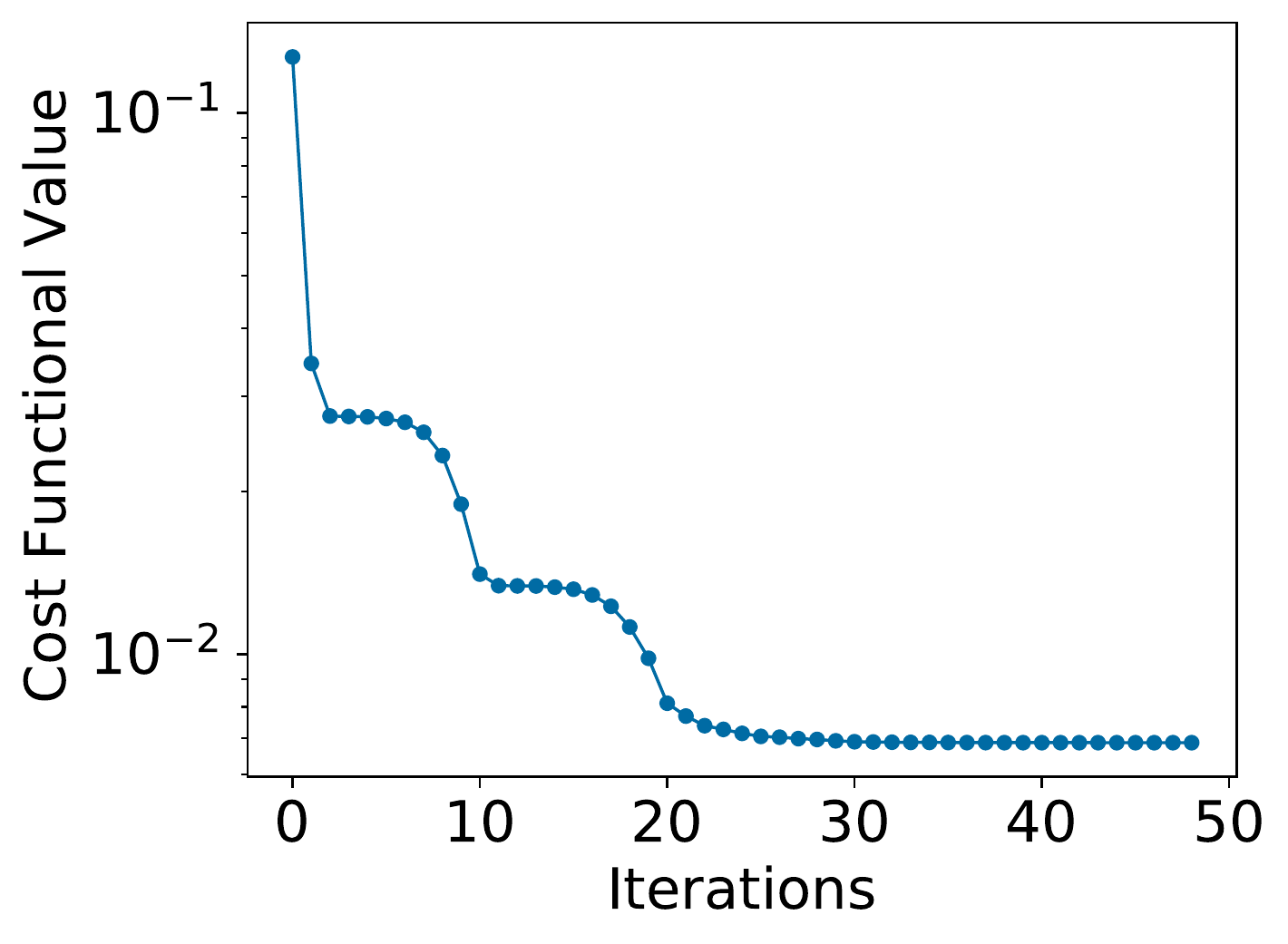}
		\caption{Cost functional value.}
		\label{sfig:id_cost_functional}
	\end{subfigure}
	\hfil
	\begin{subfigure}[t]{0.49\textwidth}
		\centering
		\includegraphics[width=0.9\textwidth]{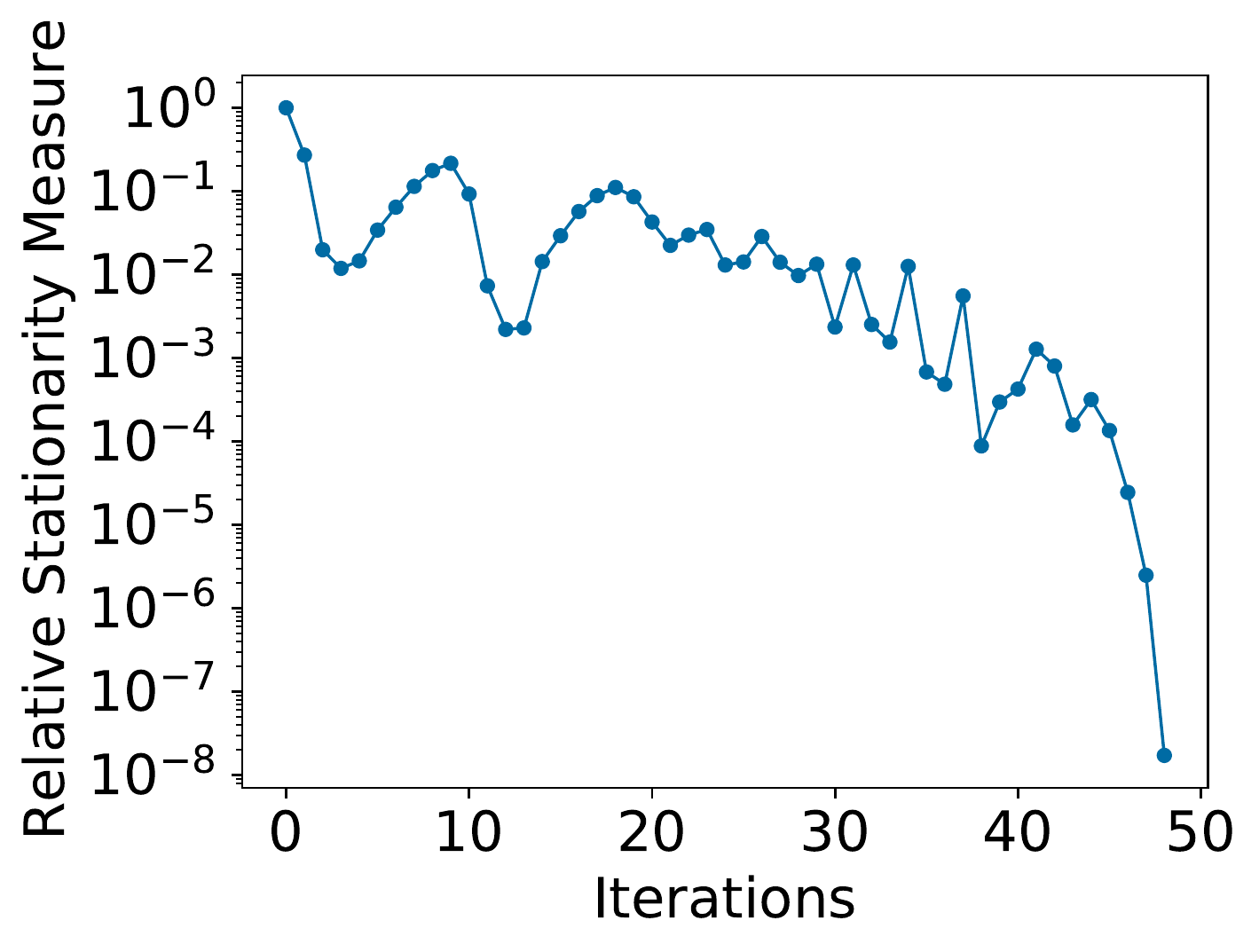}
		\caption{Relative value of the stationarity measure.}
		\label{sfig:id_stationarity_measure}
	\end{subfigure}
	\caption{History of the numerical optimization.}
	\label{fig:identification_plot}
\end{figure}

We solve the parameter identification problem \eqref{eq:parameter_identification_problem} numerically using our software cashocs described in Section~\ref{sec:adpack}. For the discretization of the adjoint system we choose analogous finite elements as for the discretization of the state system (cf. Section~\ref{ssec:numerics}). Moreover, we use a projected limited memory BFGS (L-BFGS) method with a memory of \num{5} vectors as optimization algorithm. For the computation of the step size we use an Armijo line search with initial guess one for the step size. Moreover, we restart the algorithm with a gradient step in case the respective curvature condition for the BFGS method is not satisfied. The optimization algorithm is terminated once a relative tolerance of \num{1e-6} for the stationarity measure is reached. Our initial guess for the kinetic parameters is based on the parameters computed in \cite{Engelbrecht2017Carbon, Engelbrecht2017Experimentation} and is given by $\ea = \SI{65}{\kilo\joule \per \mole}$, $\logA = 12$, and $\nexp = \num{0.222}$. We refer the reader, e.g., to \cite{Kelley1999Iterative, Nocedal2006Numerical} for a detailed description of the algorithm and to \cite{Andres2020Identification} for its application in the context of PDE constrained optimization.

\begin{figure}[!t]
	\centering
	\includegraphics[width=0.8\textwidth]{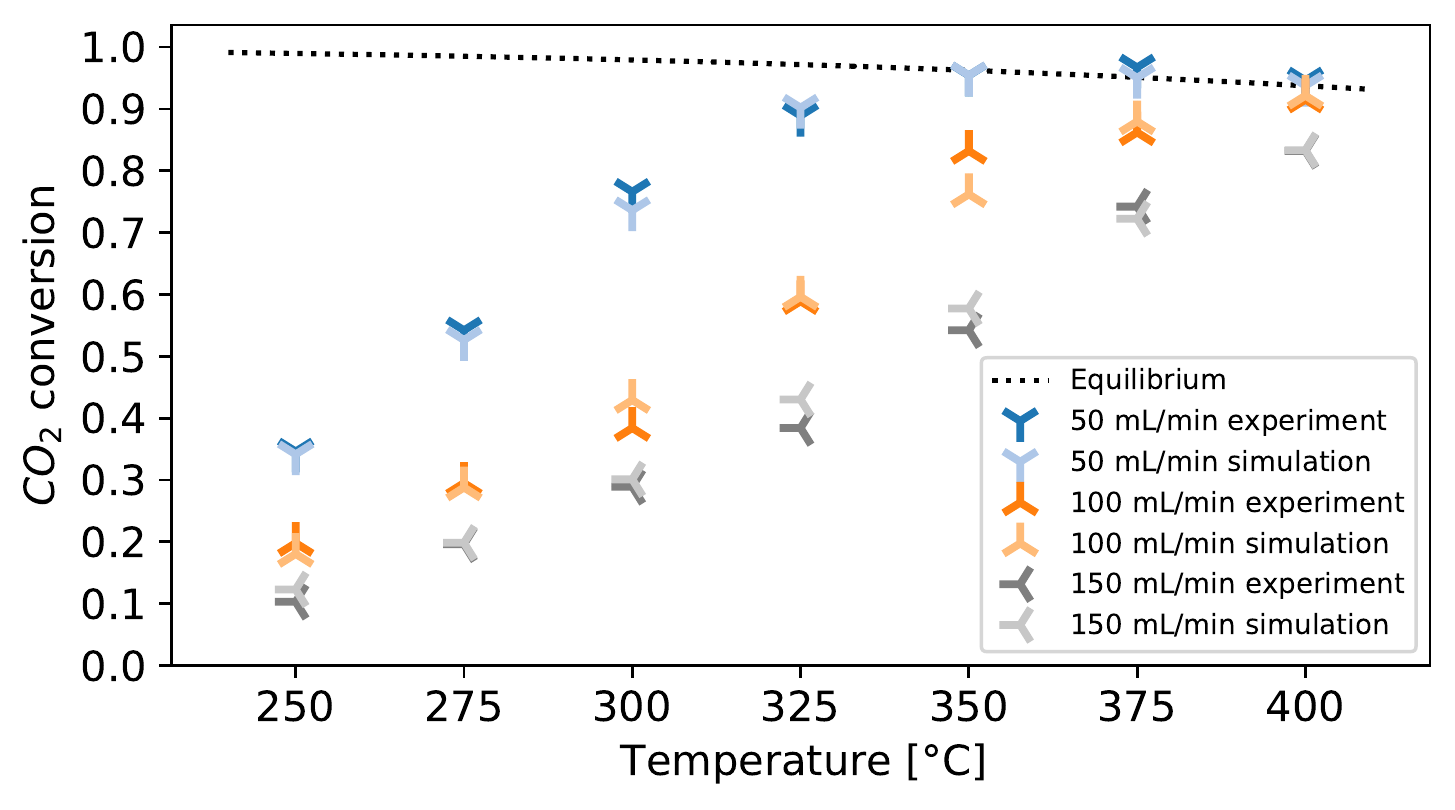}
	\caption{Simulated and measured \ce{CO2} conversion for the experiments from \cite{Engelbrecht2017Carbon, Engelbrecht2017Experimentation}.}
	\label{fig:identification}
\end{figure}

The history of the cost functional and the relative stationarity measure over the course of the optimization are shown in Figure~\ref{fig:identification_plot}. We observe that the algorithm terminates successfully after 48 iterations, suggesting that we find a local minimizer (or stationary point) of \eqref{eq:parameter_identification_problem}. Since one solve of the state system requires 21 solves of \eqref{eq:pde1d}, one for each experiment, this amounts to a total of 1113 solves for \eqref{eq:pde1d}, and 1008 solves for the corresponding adjoint system over the course of the optimization. Note, that the additional solves needed for the state system are a consequence of the Armijo line search.
%The L-BFGS method terminated after 48 iterations with the optimized values presented in Table~\ref{tab:identified_parameters}. First, we observe from the history of the cost functional and stationarity measure plotted in Figure~\ref{fig:identification_plot} that the optimization was terminated successfully, having at least found a stationary point of problem \eqref{eq:parameter_identification_problem} since the cost functional decreased by nearly a factor of 20 and the stationarity measure reached the relative tolerance of \num{1e-6}. 
In Figure~\ref{fig:identification}, the simulated and measured \ce{CO2} conversions are depicted for all 21 experimental settings for the identified kinetic parameters. We see that our model shows excellent agreement to the experimental results. The largest difference between simulated and measured \ce{CO2} conversion is about \num{7}\ \%, and the mean error is around \num{1.9}\ \%, which is well within the reproducibility of the experiment of \num{5}\ \% from \cite{Engelbrecht2017Experimentation, Engelbrecht2017Carbon} and closer than the original fit proposed there.
%Additionally, the results we obtained with our model seem to fit the experiments even more closely than the original simulation results of \cite{Engelbrecht2017Experimentation, Engelbrecht2017Carbon}, particularly, in that they behave correctly w.r.t.\the thermodynamic equilibrium of the reaction, which is not always the case for the original results. 
For these reasons, we conclude that our model is able to simulate the physical and chemical processes in the reactor sufficiently well so that we can use it in Section~\ref{sec:optimal_control_reaction} to optimize the reactor.

\begin{table}[!b]
	\centering
	{\footnotesize
		\rowcolors{2}{\tablegray}{white}
		\setlength{\tabcolsep}{1em}
		\begin{tabular}{l S}
			\toprule
			parameter [unit] & {value} \\
			\midrule
			empirical exponent $\nexp$ & 0.0581\\
			pre-exponential factor $\preexp$ [$\nicefrac{1}{\si{\second}} (\nicefrac{\si{\mole}}{\si{\cubic \meter}})^{1 - 5\nexp}$] & 4.744e5 \\
			activation energy $\ea$ [\si{\kilo\joule\per\mole}] & 52.141 \\
			\bottomrule
		\end{tabular}
		\caption{Kinetic parameters determined by the parameter identification.}
		\label{tab:identified_parameters}
	}
\end{table}

We also remark that for the parameters identified above, the solution shows a nearly isothermal behavior, in particular, the difference between the simulated gas temperature and the wall temperature is well below \SI{1}{\celsius} at all points in the reactor. Hence, we could also use an isothermal model for the reactor, where the temperature is not treated as a state variable, but directly prescribed as the wall temperature, simplifying the numerical simulations. 
%This is, partially, due to the fact that the value of the empirical exponent $\nexp$ determined through \eqref{eq:parameter_identification_problem} is comparatively small so that the reaction is distributed along the channel and we do not have temperature peaks at the beginning of the channel. 
However, as this behavior cannot be known before the determination of the kinetic parameters and since other sets of parameters could yield a non-isothermal behavior, we keep the temperature as a state variable.
%This additionally highlights the capabilities of our software package to deal with highly nonlinear and coupled PDE systems. 

Finally, the identified parameters are shown in Table~\ref{tab:identified_parameters}. We remark that they are in a good agreement with the kinetic parameters found in \cite{Falbo2018Kinetics}, where the empirical exponent $\nexp$ and the activation energy $\ea$ were found to be \num{0.076} and \SI{65.2}{\kilo\joule\per\mole}. Since we consider a higher operating pressure, a different type of catalyst, and also use a slightly different model to the ones in \cite{Falbo2018Kinetics}, the minor deviations for the identified parameters are justified.

%We remark that they deviate slightly from values of the kinetic parameters found in the literature (see, e.g., \cite{Falbo2018Kinetics} for an overview). However, this can be explained by the differences between our model and the ones considered there. In particular, we use a homogenized model for the chemical reaction and also simulate an inlet stage $\Omega\subin$, where no reaction occurs, which none of the models in \cite{Falbo2018Kinetics} do, so that the deviations for the obtained parameters are justified.

\subsection{Numerical Comparison of the Reactor Models}
\label{ssec:comparison}

To conclude this section, we briefly compare the one-dimensional model \eqref{eq:pde1d} to the three-dimensional model \eqref{eq:pde_system}. For this, we simulate all of the 21 test cases of \cite{Engelbrecht2017Experimentation, Engelbrecht2017Carbon}, using the kinetic parameters obtained previously, with both models and compute the corresponding relative errors between the models in the $L^\infty(\Omega), L^2(\Omega)$ and $L^1(\Omega)$ norms.
%\begin{table}[!b]
%	\centering
%	{\footnotesize
%		\rowcolors{2}{\tablegray}{white}
%		\setlength{\tabcolsep}{1em}
%		\begin{tabular}{c c c c c c c c}
%			\toprule
%			& $\velocity$ & $\pressure$ & $\temperature$ & $Y_{\ce{CO2}}$ & $Y_{\ce{H2}}$ & $Y_{\ce{CH4}}$ & $Y_{\ce{H2O}}$ \\
%			\midrule
%			$L^\infty$ mean & \num{0.069}\ \% & \num{0.200}\ \% & \num{0.002}\ \% & \num{0.018}\ \% & \num{0.018}\ \% & \num{0.032}\ \% & \num{0.032}\ \% \\
%			%
%			$L^\infty$ max  & \num{0.079}\ \% & \num{0.209}\ \% & \num{0.004}\ \% & \num{0.048}\ \% & \num{0.052}\ \% & \num{0.053}\ \% & \num{0.052}\ \% \\
%			\midrule
%			$L^2$ mean & \num{0.066}\ \% & \num{0.201}\ \% & \num{0.001}\ \% & \num{0.020}\ \% & \num{0.018}\ \% & \num{0.036}\ \% & \num{0.036}\ \% \\
%			%
%			$L^2$ max  & \num{0.070}\ \% & \num{0.21}\ \% & \num{0.003}\ \% & \num{0.043}\ \% & \num{0.042}\ \% & \num{0.053}\ \% & \num{0.053}\ \% \\
%			\midrule
%			$L^1$ mean & \num{0.066}\ \% & \num{0.201}\ \% & \num{0.001}\ \% & \num{0.020}\ \% & \num{0.018}\ \% & \num{0.036}\ \% & \num{0.036}\ \% \\
%			%
%			$L^1$ max  & \num{0.070}\ \% & \num{0.210}\ \% & \num{0.003}\ \% & \num{0.043}\ \% & \num{0.042}\ \% & \num{0.053}\ \% & \num{0.053}\ \% \\
%			\bottomrule
%		\end{tabular}
%		\caption{Average and maximum relative errors between the one-dimensional porous medium model and the three-dimensional model different norms.}
%		\label{tab:error_analysis}
%	}
%\end{table}
\begin{figure}[!t]
	\centering
	\includegraphics[width=0.75\textwidth]{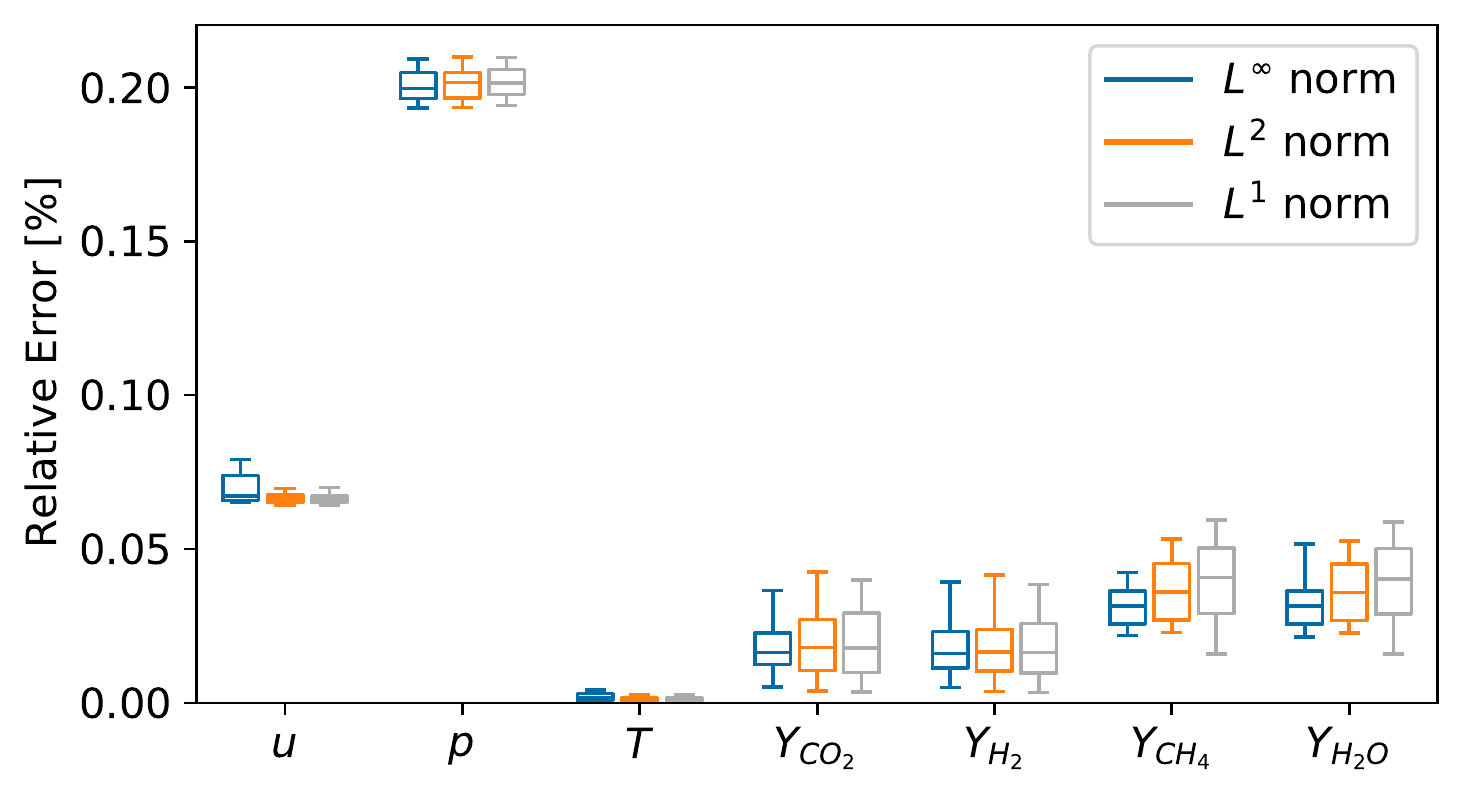}
	\caption{Relative errors between the 1D and 3D model.}
	\label{fig:error_analysis}
\end{figure}
To compare the models to each other, we average the results of the three-dimensional model over the channel cross section, which can then be directly compared to the results of the one-dimensional model. The arising errors are shown in Figure~\ref{fig:error_analysis} in the form of a box plot. We see that the relative error between both models is well below \num{1}\ \% for all test cases and variables. The highest relative errors are obtained for the pressure with about \num{0.2}\ \%, all other variables have relative errors even below 1~\textperthousand. This is the case since the physical and chemical behavior in a single channel of the reactor is basically one-dimensional, with variations along the channel length being significantly larger than variations along its cross section. Finally, we note that the relative errors for the temperature are very small even compared to the small errors of the other variables. This is due to fact that we use the Kelvin scale to compute the relative errors for the temperature as this is a ratio scale. These results validate our approach to use the 1D porous medium model for computing the necessary kinetic parameters for the models instead of the more costly three-dimensional one since both yield basically identical results.

\section{Optimization of the Microreactor}
\label{sec:optimal_control_reaction}

In this section, we discuss the optimization of the microchannel reactor under consideration by means of solving PDE constrained optimal control problems. Throughout this section, we use the one-dimensional model \eqref{eq:pde1d}, which models the behavior of the reactor as validated in the previous section, as corresponding state system for the optimization problems. We consider two optimization problems for improving the reactor's performance which we solve numerically using our software package cashocs described in Section~\ref{sec:adpack}. We discuss the obtained results as well as their applicability and realizability.

\subsection{Optimization with Fixed Flow Rates}
\label{ssec:fixed_flow_opt}

For the first problem, we keep the inlet flow rate fixed at the respective values considered in \cite{Engelbrecht2017Experimentation, Engelbrecht2017Carbon}. Our goal is to maximize the \ce{CO2} conversion of the reactor, since this corresponds to maximizing the product yield for a fixed flow rate, by using the wall temperature $\temperature\subwall$ as optimization variable, which is detailed below. To model this, we consider the following optimization problem
\begin{equation}
	\label{eq:opt_cc}
	\min_{y, \temperature\subwall}\ J(y, \temperature\subwall) = \frac{1}{2} \integral{\Gamma\subout} \left( \conv{\ce{CO2}} - 1 \right)^2 \dmeas{s} \quad \text{ s.t. } \quad e(y, \temperature\subwall) = 0 \quad \text{ and } \quad \temperature\subwall \in U_\text{ad}.
\end{equation}
As in Section~\ref{sec:parameter_identification}, the state variables are combined into the vector $y = [\pressure, \velocity, \temperature, \massfracvec]\transposed$ and the \ce{CO2} conversion $\conv{\ce{CO2}}$ is defined as in \eqref{eq:def_conversion}. Moreover, the wall temperature $\temperature\subwall$ now plays the role of the control variable, and the state system $e(y, \temperature\subwall) = 0$ is given by the corresponding weak form of \eqref{eq:pde1d}, with the respective inlet conditions for the velocity. Note, that we use a tracking-type cost functional which aims at bringing the \ce{CO2} conversion as close as possible to \num{1} on $\Gamma\subout$. Note, that this case, i.e., $\conv{\ce{CO2}} = 1$, corresponds to the total conversion of \ce{CO2} and, therefore, to the maximum product yield. We remark that we tested several alternative cost functions that also have the goal of maximizing the reactor's yield and found that all of them gave nearly identical results in practice. Hence, it is justified that we restrict ourselves to the cost functional given in \eqref{eq:opt_cc}. Finally, the set of admissible controls $U_\text{ad}$ is given by
\begin{equation*}
	U_\text{ad} = \Set{\temperature\subwall \in L^2(\Omega) | \SI{180}{\celsius} \leq \temperature\subwall \leq \SI{600}{\celsius} \text{ a.e. in } \Omega},
\end{equation*}
i.e., we have box constraints for the wall temperature. Due to these box constraints, we do not require additional regularization terms for the optimization problem \eqref{eq:opt_cc}. Since the model is quasi-isothermal (cf. Section~\ref{ssec:results_identification}), this constraint also restricts the temperature of the gas mixture. Note, that the lower bound of \SI{180}{\celsius} corresponds to the boiling point of water at a pressure of \SI{10}{\bar}, so that we do not have to consider any phase change effects. The upper bound for the temperature is chosen so that we can neglect the effect of the RWGS reaction and the corresponding \ce{CO} production (cf. Section~\ref{ssec:sabatier_reaction}).

%\begin{table}[!b]
%	\centering
%	{\footnotesize
%		\rowcolors{2}{\tablegray}{white}
%		\setlength{\tabcolsep}{1em}
%		\begin{tabular}{r c c c c}
%			\toprule
%			inlet flow rate & constant & two stages & three stages & distributed \\
%			\midrule
%			\SI{50}{\milli\liter \per \minute} & 6 / 6 & \phantom{0}8 / \phantom{0}8 & 12 / 12 & 13 / 13 \\
%			\SI{100}{\milli\liter \per \minute} & 9 / 9 & 10 / 10 & 12 / 12 & 18 / 18 \\
%			\SI{150}{\milli\liter \per \minute} & 8 / 8 & \phantom{0}8 / \phantom{0}8  & 10 / 10 & 16 / 15 \\
%			\bottomrule
%		\end{tabular}
%		\caption{Number of solves for state / adjoint system needed for the optimization.}
%		\label{tab:opt_cc_iterations}
%	}
%\end{table}

For the wall temperature, we consider the following four models. First, we use a constant value for $\temperature\subwall$ throughout the entire reactor, which corresponds to an isothermal reactor. Second, we divide $\Omega\subreac$ into two halves along its length and restrict $\temperature\subwall$ to be constant in each half, which models a reactor consisting of two isothermal stages. Third, we proceed analogously to before, but now divide $\Omega\subreac$ into thirds, which corresponds to a three-stage isothermal reactor. For the final model, we do not specify any sort of profile and consider $\temperature\subwall$ to be an arbitrary function in $L^2(\Omega)$, which we discretize with piecewise linear Lagrangian finite elements. We refer to these models as constant, two stage, three stage, and distributed model, respectively, throughout the rest of this paper. These models exhibit different levels of complexity, particularly regarding their realizability and capability of controlling the reaction, which we discuss in Section~\ref{ssec:realizability}. In particular, we note that while these models may only be approximately realizable in practice, they still serve as benchmark for the optimization by showcasing the theoretical limitations of the reactor.

\begin{table}[!b]
	\centering
	{\footnotesize
		\rowcolors{2}{\tablegray}{white}
		\setlength{\tabcolsep}{1em}
		\begin{tabular}{r c c c c}
			\toprule
			inlet flow rate & constant & two stages & three stages & distributed \\
			\midrule
			\SI{50}{\milli\liter \per \minute} & 6 / 6 & \phantom{0}8 / \phantom{0}8 & 12 / 12 & 13 / 13 \\
			\SI{100}{\milli\liter \per \minute} & 9 / 9 & 10 / 10 & 12 / 12 & 18 / 18 \\
			\SI{150}{\milli\liter \per \minute} & 8 / 8 & \phantom{0}8 / \phantom{0}8  & 10 / 10 & 16 / 15 \\
			\bottomrule
		\end{tabular}
		\caption{Number of solves for state / adjoint system for problem \eqref{eq:opt_cc}.}
		\label{tab:opt_cc_iterations}
	}
\end{table}

\begin{figure}[!b]
	\centering
	\begin{subfigure}{0.49\textwidth}
		\centering
		\includegraphics[width=0.9\textwidth]{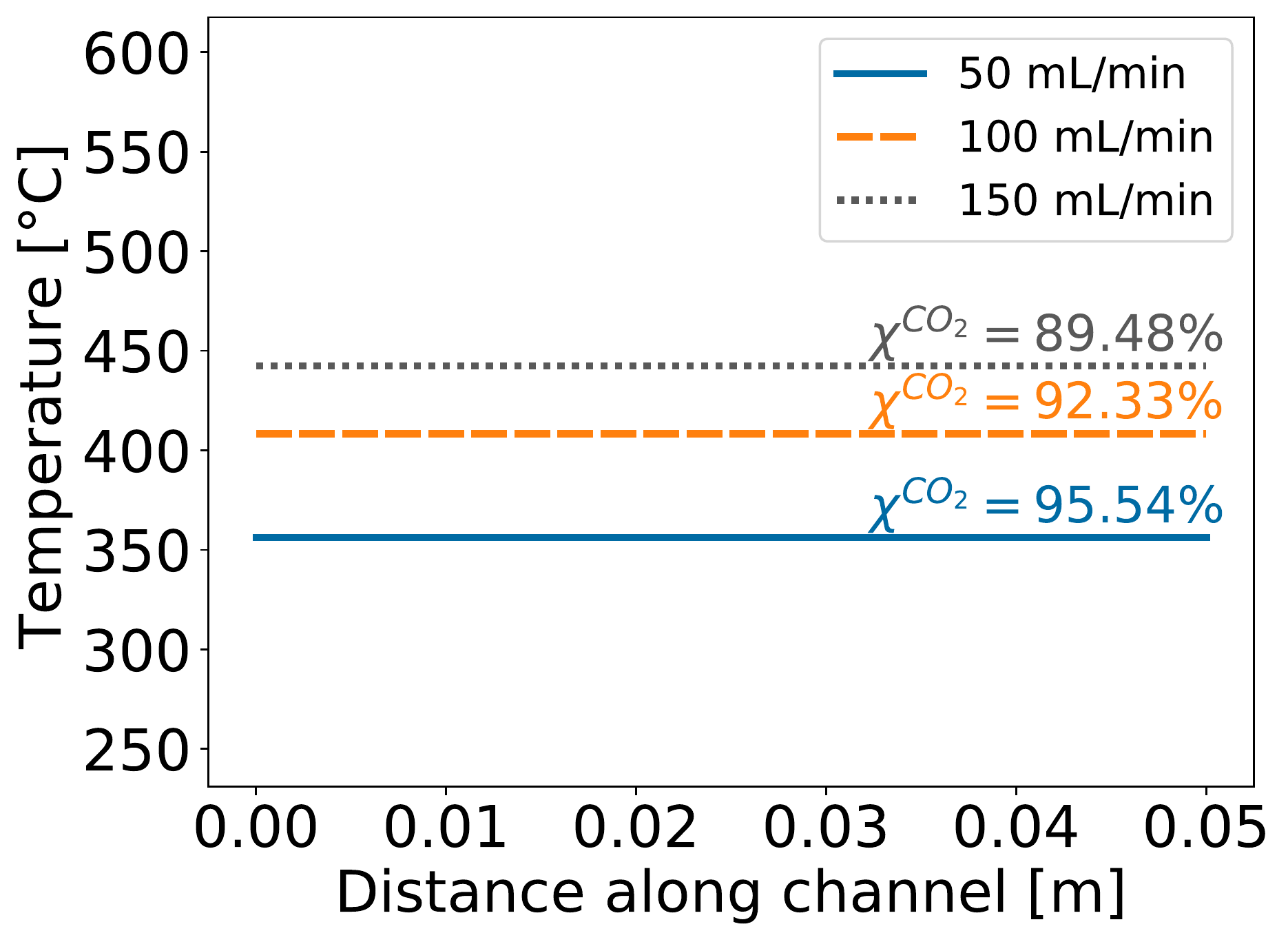}
		\caption{Constant temperature model.}
	\end{subfigure}
	\hfil
	\begin{subfigure}{0.49\textwidth}
		\centering
		\includegraphics[width=0.9\textwidth]{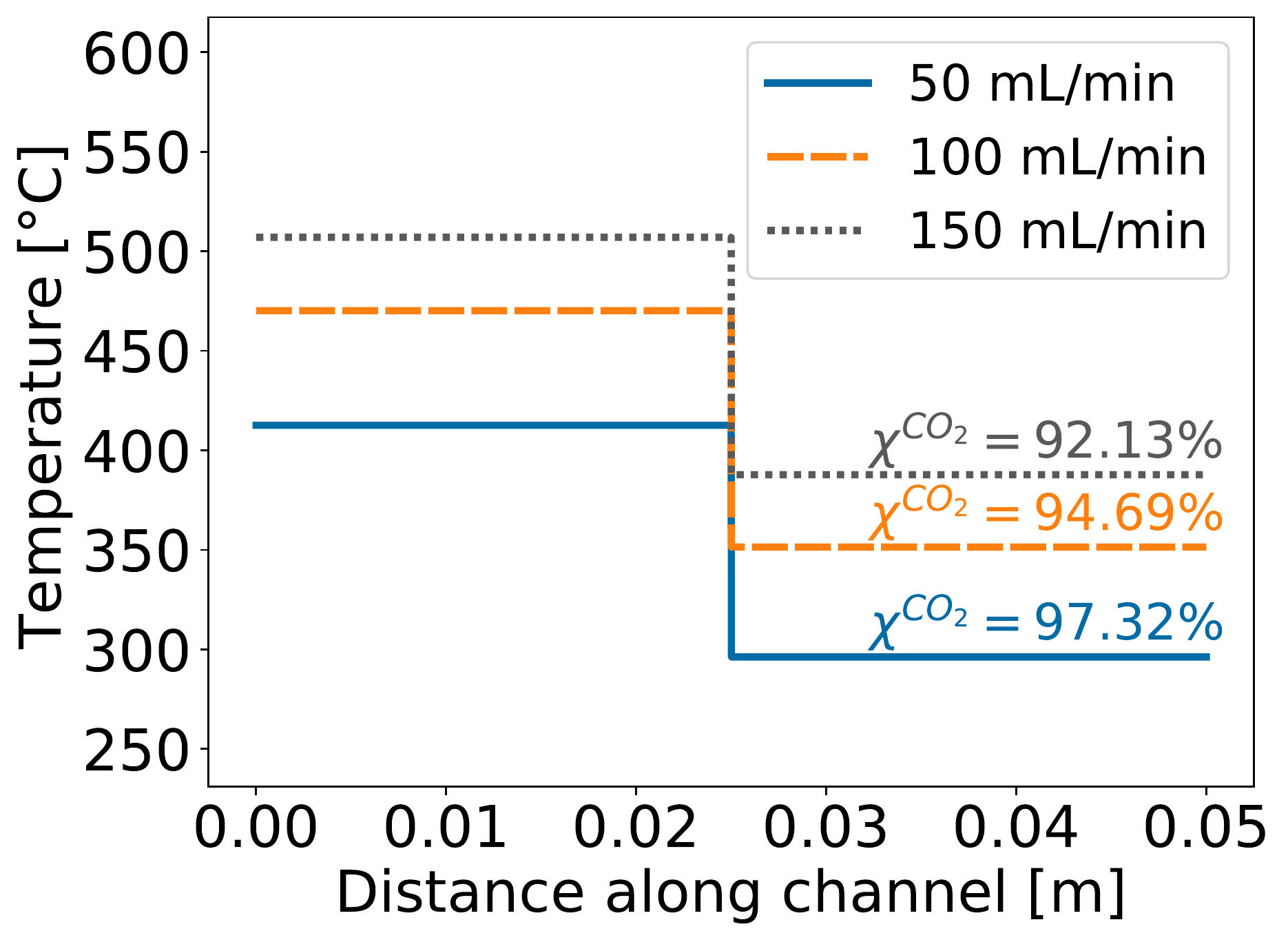}
		\caption{Two temperature stages model.}
	\end{subfigure}
	\vskip\baselineskip
	\begin{subfigure}{0.49\textwidth}
		\centering
		\includegraphics[width=0.9\textwidth]{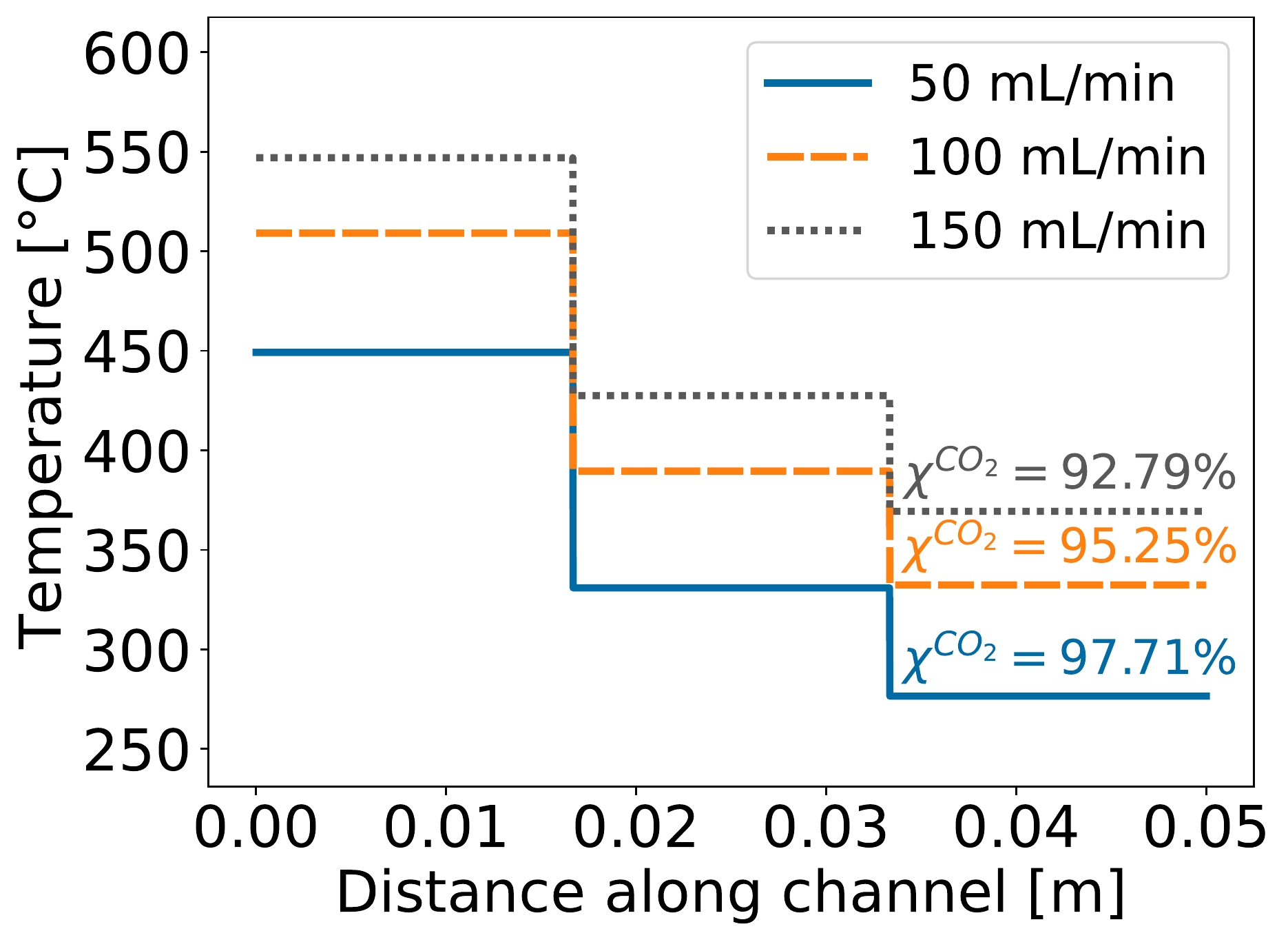}
		\caption{Three temperature stages model.}
	\end{subfigure}
	\hfil
	\begin{subfigure}{0.49\textwidth}
		\centering
		\includegraphics[width=0.9\textwidth]{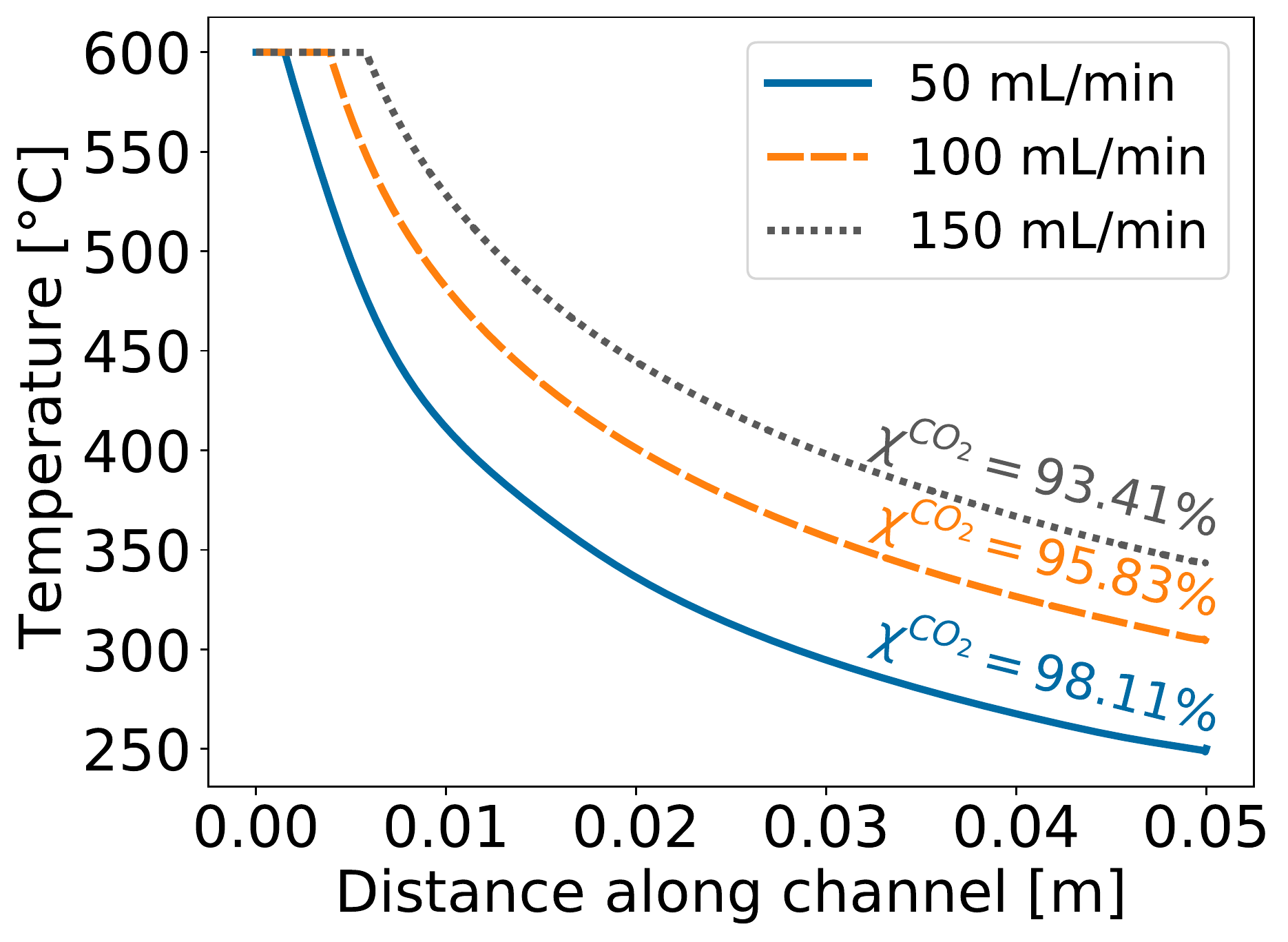}
		\caption{Distributed temperature model.}
	\end{subfigure}
	\caption{Optimized temperature profiles for problem \eqref{eq:opt_cc}.}
	\label{fig:opt_cc}
\end{figure}

We solve the optimization problem \eqref{eq:opt_cc} numerically for all four models and the three inlet flow rates of \SI{50}{\milli\liter\per \minute}, \SI{100}{\milli\liter\per \minute}, and \SI{150}{\milli\liter\per \minute}, as considered in \cite{Engelbrecht2017Experimentation,Engelbrecht2017Carbon}, using our software package cashocs (cf. Section~\ref{sec:adpack}). For the corresponding adjoint system we use an analogous discretization to the one for the state system (cf. Section~\ref{ssec:numerics}). As optimization algorithm we again use the projected L-BFGS method described in Section~\ref{ssec:results_identification}, which is terminated once the relative stationarity measure reaches a tolerance of \num{1e-4}. The required amount of PDE solves for the state and adjoint systems for the optimization algorithm are tabulated in Table~\ref{tab:opt_cc_iterations}. Note, that the number of iterations of the algorithm corresponds to the number of solves for the adjoint system minus one, as we need an additional gradient computation for the termination criterion. We observe that the algorithm is very efficient as it needs less than 20 iterations to reach the desired tolerance for the stationarity measure in all cases. Since the number of solves for the state system exceeds the ones for the adjoint system only for one case and one additional solve, we observe that the initial step length is almost always accepted by the algorithm. We note that with increasing complexity of the models we need slightly more iterations to reach the desired tolerance which is, however, to be expected.

The optimized wall temperatures are depicted in Figure~\ref{fig:opt_cc} together with the corresponding \ce{CO2} conversions. Comparing the \ce{CO2} conversions with the ones from the experiments (cf. Section~\ref{sec:parameter_identification}) we recognize that we achieve considerably higher levels of \ce{CO2} conversion for all flow rates by optimizing the wall temperature. Moreover, by comparing the achieved \ce{CO2} conversions with the equilibrium conversions (cf. Figure~\ref{sfig:sab_eq_conv}), we can also conclude that the thermodynamic equilibrium is a limiting factor in all cases as the difference between the obtained and the equilibrium \ce{CO2} conversions are only about \num{1}--\num{2}\ \% for the respective outlet temperatures. We observe that the optimized wall temperatures are at their largest near the inlet of the channel and then decrease monotonically along its length. This yields a large reaction rate at the beginning of the reactor as well as more favorable thermodynamic equilibria towards its end, which are actually limiting the reaction as discussed above. Moreover, note that the optimized wall temperatures increase monotonically with the flow rate, due to the following. As a higher flow rate corresponds to a lower residence time of the gas in the reactor, a higher reaction rate is needed to increase the \ce{CO2} conversion, which can be achieved through a higher wall temperature (cf. Figure~\ref{sfig:sab_reac_rates}). This confirms our deliberations from Section~\ref{ssec:sabatier_reaction} regarding the dependence of the thermodynamic equilibrium and the reaction rate on the temperature w.r.t.\ the optimization of the reactor.

Analyzing the obtained \ce{CO2} conversions more closely, we first observe that the \ce{CO2} conversion increases monotonically with the complexity of the wall temperature models, i.e., the constant model performs worst for all flow rates, the two and three stage models perform increasingly better, and the distributed model has the highest conversion of all profile types for each flow rate. In particular, the increase in conversion is roughly \num{2}\ \% between the constant and two level profiles and then only about \num{0.5}\ \% each for the successive models. Hence, we conclude that the ability to use a lower temperature towards the end of the reactor can have an important effect on the quality of the reaction, but using more sophisticated temperature profiles does not add a great amount of additional improvements on top of this. Moreover, we observe that the optimized conversion also decreases monotonically with increasing flow rate. This is due to the fact that for higher flow rates a larger temperature has to be used to obtain a sufficiently high reaction rate, as mentioned previously. This then leads to worse thermodynamic equilibria and explains the overall lower \ce{CO2} conversion. 

\begin{table}[!t]
	\centering
	{\footnotesize
		\rowcolors{2}{\tablegray}{white}
		\setlength{\tabcolsep}{1em}
		\begin{tabular}{c c c c c}
			\toprule
			inlet flow rate & constant & two stages & three stages & distributed \\
			\midrule
			\SI{50}{\milli\liter \per \minute} & \num{2.131e-05} & \num{2.171e-05} & \num{2.180e-05} & \num{2.189e-05} \\
			\SI{100}{\milli\liter \per \minute} & \num{4.119e-05} & \num{4.225e-05} & \num{4.250e-05} & \num{4.275e-05} \\
			\SI{150}{\milli\liter \per \minute} & \num{5.988e-05} & \num{6.166e-05} & \num{6.210e-05} & \num{6.251e-05} \\
			\bottomrule
		\end{tabular}	
	}
	\caption{Molar product yield in [\si{\mol \per \second}] for problem \eqref{eq:opt_cc}.}
	\label{tab:opt_cc_yield}
\end{table}

Finally, let us investigate the obtained product yield which is shown in Table~\ref{tab:opt_cc_yield}. There, we consider the total amount of product, i.e., \ce{CH4} and \ce{H2O}, that is generated by the reactor and, due to the stoichiometry of the reaction, \nicefrac{1}{3} of this amount is \ce{CH4} and \nicefrac{2}{3} is \ce{H2O}. Obviously, for a fixed flow rate a higher conversion leads to a higher molar product yield, so that the yield increases monotonically with the model complexity. Moreover, for the conditions under consideration, we observe that the molar yield is almost linear in the flow rate for all models. This can be explained by the fact that the achieved conversions are all rather close to each other for all models and flow rates. Therefore, the slightly lower \ce{CO2} conversion obtained at higher flow rates is compensated by the higher resulting throughput, leading to overall higher product yields for higher flow rates. These results indicate that one should try to use the highest possible flow rate so that the resulting quality of the product, measured by the \ce{CO2} conversion, is still sufficiently high. We take a more detailed look at this problem of balancing the molar yield and product quality of the reactor in the following section.

%Altogether, we observe that optimizing the wall temperature for a fixed flow rate can significantly increase the performance of the reactor. In particular, we have shown that already optimizing a constant wall temperature can improve the obtained \ce{CO2} conversion significantly and yield better reactor performance. Additionally, using a two-stage temperature model that allows the temperature to decrease at the end of the reactor leads to better thermodynamic equilibria and, thus, further increases the \ce{CO2} conversion. Using even more sophisticated temperature models can enhance the reactors performance, but the corresponding temperature profiles can also be harder to realize in practice.

\subsection{Optimization with Variable Flow Rates and State Constraints}
\label{ssec:state_constraint}

As indicated above, our goal is now to maximize the molar yield of the reaction by considering the wall temperature and the inlet flow rate as control variables, where the product quality is ensured by means of a state constraint on the \ce{CO2} conversion. In particular, we consider the following optimization problem
\begin{equation}
	\label{eq:opt_sc_theo}
	\min_{y, u_c}\ \costfunction(y, u_c) = -\integral{\Gamma\subout} \density \velocity \cdot \normal \dmeas{s} \quad \text{ s.t. } \quad e(y, u_c) = 0, \quad y\in Y_\text{ad}, \quad \text{ and } \quad u_c \in U_\text{ad}.
\end{equation}
Analogously to before, the vector of state variables is given by $y = [\pressure, \velocity, \temperature, \massfracvec]\transposed$ and the vector of control variables is given by $u_c = [\velocity\subin, \temperature\subwall]\transposed$, i.e., now both the inlet flow rate and the wall temperature act as control variables. The state system $e(y, u_c) = 0$ is, again, given by the corresponding weak form of the one-dimensional porous medium model \eqref{eq:pde1d}. The cost functional $J$ measures the negative outgoing mass flow rate at the outlet of the reactor. Note, that thanks to the sign in $J$, the minimization of this cost functional is equivalent to maximizing the outgoing mass flow rate. The set of admissible controls is defined as
\begin{equation*}
	U_\text{ad} = \Set{[\velocity\subin, \temperature\subwall]\transposed \in \R \times L^2(\Omega) | u_a \leq \velocity\subin \leq u_b \text{ on } \Gamma\subin \quad \text{ and } \quad \SI{180}{\celsius} \leq \temperature\subwall \leq \SI{600}{\celsius} \text{ a.e. in } \Omega},
\end{equation*}
where the box constraints for $\temperature\subwall$ are the same as in the previous section. The constraints for the inlet flow rate are motivated by the physical limitations of the reactor, i.e., we cannot prescribe arbitrarily small or large flow rates as the reactor is not stable for such conditions. Note, that the particular choice of $u_a$ and $u_b$ is detailed below. Additionally, as we use the one-dimensional model of the reactor, the inlet velocity $\velocity\subin$ corresponds to the averaged inlet flow rate, and is, thus, a scalar quantity. The state constraint is modeled by the set of admissible states $Y_\text{ad}$ which is given by
\begin{equation*}
	Y_\text{ad} = \Set{[\pressure, \velocity, \temperature, \massfracvec]\transposed | \conv{\ce{CO2}} = 1 - \frac{Y_{\ce{CO2}}}{Y\subin_{\ce{CO2}}} \geq \scconst \text{ on } \Gamma\subout},
\end{equation*}
i.e., the state constraint is given by the requirement, that the \ce{CO2} conversion is larger than a given desired conversion $0 \leq \scconst < 1$. Obviously, this can be easily transformed into a more classical form of a state constraint, namely $Y_{\ce{CO2}} \leq C$ for some appropriate constant $C \in (0,1)$. As we have seen in Section~\ref{ssec:fixed_flow_opt}, the maximum achievable \ce{CO2} conversion decreases monotonically with the flow rate. Therefore, we conclude that for the maximum possible flow rate the state constraint is satisfied with equality. Since the molar product yield of the reactor solely depends on the mass flow rate $\density \velocity$ and the gas composition, which is restricted to the desired conversion as just discussed, the solution of problem \eqref{eq:opt_sc_theo} maximizes the product yield of the reactor for a given \ce{CO2} conversion $\scconst$.
%Additionally, we remark that we consider the mass flow rate $\density\velocity$ and not the gas velocity $\velocity$ alone since the former is independent of the temperature at the reactor outlet. 

Finally, we note that problem \eqref{eq:opt_sc_theo} is a Dirichlet boundary control problem for the control variable $\velocity\subin$. For the numerical solution of this problem we proceed analogously to \cite[Chapter 10.6]{Delfour2011Shapes}, i.e., we treat the Dirichlet boundary condition with the help of an appropriately chosen Lagrange multiplier that consists of terms involving the respective adjoint variables. These terms correspond to the ones arising from the integration by parts carried out for the derivation of the variational formulation of the PDE. This allows us to include the boundary condition in the variational form of the problem, so that we can solve it with our software cashocs. For more details regarding Dirichlet boundary control problems we refer the reader, e.g., to \cite{Hinze2009Optimization}.

For the numerical solution of problem \eqref{eq:opt_sc_theo} we treat the state constraint with a Moreau-Yosida regularization \cite{Hinze2009Optimization} which leads to the following regularized optimization problem
\begin{equation}
	\label{eq:opt_sc}
	\begin{aligned}
		\min_{u_c}\ J_\gamma(y, u_c) &= -\integral{\Gamma\subout} \density \velocity \cdot \normal \dmeas{s} + \frac{1}{2\gamma} \integral{\Gamma\subout} \max\left( 0, \hat{\mu} + \gamma \left( \scconst - \conv{\ce{CO2}} \right) \right)^2 \dmeas{s} \\
		&\text{ s.t. } \quad e(y, u_c) = 0 \quad \text{ and } \quad u_c \in U_\text{ad}.
	\end{aligned}
\end{equation}
Here, the state and control variables as well as the state system and $U_\text{ad}$ are defined exactly as for problem \eqref{eq:opt_sc_theo}. The cost functional $J_\gamma$ now includes the Moreau-Yosida regularization of the state constraint with penalty parameter $\gamma > 0$ and shift function $\hat{\mu}$. However, we only consider the case $\hat{\mu} = 0$ in this paper. The idea of the Moreau-Yosida regularization is to solve a sequence of regularized optimization problems \eqref{eq:opt_sc} with increasing values of $\gamma$ with $\gamma\to \infty$ so that we obtain a solution of the state constrained problem \eqref{eq:opt_sc_theo} in the limit. 
%Note, that this procedure also explains the need for box constraints for the inlet flow rate in the following way. For low values of $\gamma$ the regularization term does not play a large role. Hence, the optimum inlet flow rate for the problem is to be chosen as large as possible, with the temperature not being of importance. As the value of $\gamma$ is increased, the optimum flow rate becomes smaller and smaller, such that the state constraint can be satisfied. However, the numerical solution of the system is not possible anymore for arbitrarily large flow rates and, additionally, it is very hard for optimization algorithms to again decrease a very large flow velocity to some reasonable amount. Therefore, we chose the box constraints for each problem individually such that $U_\text{ad}$ is not too large and a numerical solution of the problem is still possible.

For the wall temperature we consider the same four models as in the previous section, i.e., the constant, two and three stage, and distributed model. Again, the models might only be approximately realizable in practice, but still act as benchmark for the optimization. For the state constraint we use the following values for the desired \ce{CO2} conversion
\begin{equation*}
	\scconst \in \Set{0.85, 0.875, 0.9, 0.925, 0.95}.
\end{equation*}
Additionally, the box constraints for the velocity are always chosen with $u_a$ corresponding to \SI{50}{\milli\liter\per\minute}, i.e., the lowest flow rate considered in the experiments in \cite{Engelbrecht2017Experimentation, Engelbrecht2017Carbon}, and with $u_b$ corresponding to a flow rate of at least \SI{150}{\milli\liter\per\minute}, which is explained more detailedly below. Note, that for a flow rate of \SI{50}{\milli\liter\per\minute}, the obtained optimized temperature profiles from Section~\ref{ssec:fixed_flow_opt} all have \ce{CO2} conversions of over \num{95}\ \% (cf. Figure~\ref{fig:opt_cc}), which directly makes the corresponding controls and states feasible for problem \eqref{eq:opt_sc_theo} for all considered values of $\scconst$, so that $U_\text{ad}$ and $Y_\text{ad}$ are, in fact, non-empty and we do not minimize over the empty set.

For the numerical solution of \eqref{eq:opt_sc_theo}, we use a homotopy approach, i.e., we solve \eqref{eq:opt_sc} for one value of $\gamma$ and then use the obtained solution as initial guess for \eqref{eq:opt_sc} with the next higher value of $\gamma$. In particular, for all three piecewise constant temperature models, we use the sequence $\gamma_l = 10^l, l=0,\dots,6$, and for the distributed temperature model we use the sequence $\gamma_l = 7^l, l=0,\dots,7$, as the corresponding problem is slightly harder to solve numerically. We solve the subproblems \eqref{eq:opt_sc} by means of a projected L-BFGS method with memory size 5 and a stopping tolerance of \num{1e-2} for the relative stationarity measure, as discussed in Section~\ref{ssec:results_identification}. Again, the adjoint system is discretized analogously to the state system (cf. Section~\ref{ssec:numerics}). The value of the upper bound for the inlet velocity $u_b$ is chosen as follows. We start with a value of $u_b$ corresponding to \SI{150}{\milli\liter\per\minute} and solve the optimization problem \eqref{eq:opt_sc_theo} as discussed above. If the resulting optimized inlet velocity is located at the upper bound $u_b$, we increase the bound by \SI{50}{\milli\liter\per\minute} and repeat this procedure, until the optimized flow rate is not at the upper bound anymore. Since the maximum possible \ce{CO2} conversion decreases monotonically with the flow rate as discussed previously, it is straightforward to see that this indeed yields the maximum flow rate. Note, that we do not use a fixed and sufficiently large value for $u_b$ for all settings as we cannot solve all of the corresponding problems numerically due to the following. For the initial low values of $\gamma$ the solution to \eqref{eq:opt_sc} is given by $\velocity\subin = u_b$ as the regularization term for the state constraint has only a minor influence on the problem. The higher the value of $\gamma$, the larger this influence becomes and the smaller the flow rate has to be to satisfy the state constraint. This, however, leads to very bad initial guesses for the subsequent problem with a larger penalty parameter, making the numerical solution of \eqref{eq:opt_sc} very hard or even impossible, so that the whole approach of regularizing the state constraint becomes pointless in this case.

Before we discuss the numerical results, we remark that the absolute difference between the \ce{CO2} conversion obtained by solving \eqref{eq:opt_sc_theo} numerically and the desired conversion $\scconst$ is less than 0.01\ \% in all considered cases, which is sufficiently accurate for our purposes and shows that we indeed compute feasible solutions to the original problem \eqref{eq:opt_sc_theo}. The required amount of solves for the state and adjoint system for the numerical solution of \eqref{eq:opt_sc_theo} is shown in Table~\ref{tab:opt_sc_iterations}. Since we now solve a sequence of optimization problems, the required number of PDE solves obviously increases compared to the problems investigated in Section~\ref{ssec:fixed_flow_opt}. We also observe that we need significantly more solves for the state system than for the adjoint system which indicates that the step size computation via the Armijo rule takes more effort. On average we need about two to three times as many solves for the state system compared to the number of solves for the adjoint system. As before (cf. Table~\ref{tab:opt_cc_iterations}), we observe that we need roughly the same amount of iterations to solve problems with the piecewise constant temperature models, and that we need considerably more iterations to solve the problems for the distributed model. Finally, we note that the number of PDE solves increases with increasing $\scconst$ for almost all cases, which indicates that the problems become increasingly harder to solve for higher values of $\scconst$.  

\begin{table}[!t]
	\centering
	{\footnotesize
		\rowcolors{2}{\tablegray}{white}
		\setlength{\tabcolsep}{1em}
		\begin{tabular}{c c c c c}
			\toprule
			{\ce{CO2} conversion} & {constant} & {two stages} & {three stages} & {distributed} \\
			\midrule
			\num{85.0}\ \% & 159 / \phantom{0}71 & 147 / \phantom{0}59 & 126 / \phantom{0}48 & 247 / 152  \\
			\num{87.5}\ \% & 163 / \phantom{0}83 & 152 / \phantom{0}61 & 145 / \phantom{0}67 & 265 / 149 \\
			\num{90.0}\ \% & 166 / \phantom{0}78 & 164 / \phantom{0}67 & 146 / \phantom{0}66 & 295 / 182 \\
			\num{92.5}\ \% & 182 / \phantom{0}90 & 249 / 103 & 186 / \phantom{0}89 & 331 / 198 \\
			\num{95.0}\ \% & 215 / 124 & 201 / \phantom{0}97 & 194 / 105 & 254 / 141 \\
			\bottomrule
		\end{tabular}
		\caption{Number of solves for state / adjoint system for problem \eqref{eq:opt_sc_theo}.}
		\label{tab:opt_sc_iterations}
	}
\end{table}

\begin{figure}[!b]
	\centering
	\begin{subfigure}[t]{0.49\textwidth}
		\centering
		\includegraphics[width=0.9\textwidth]{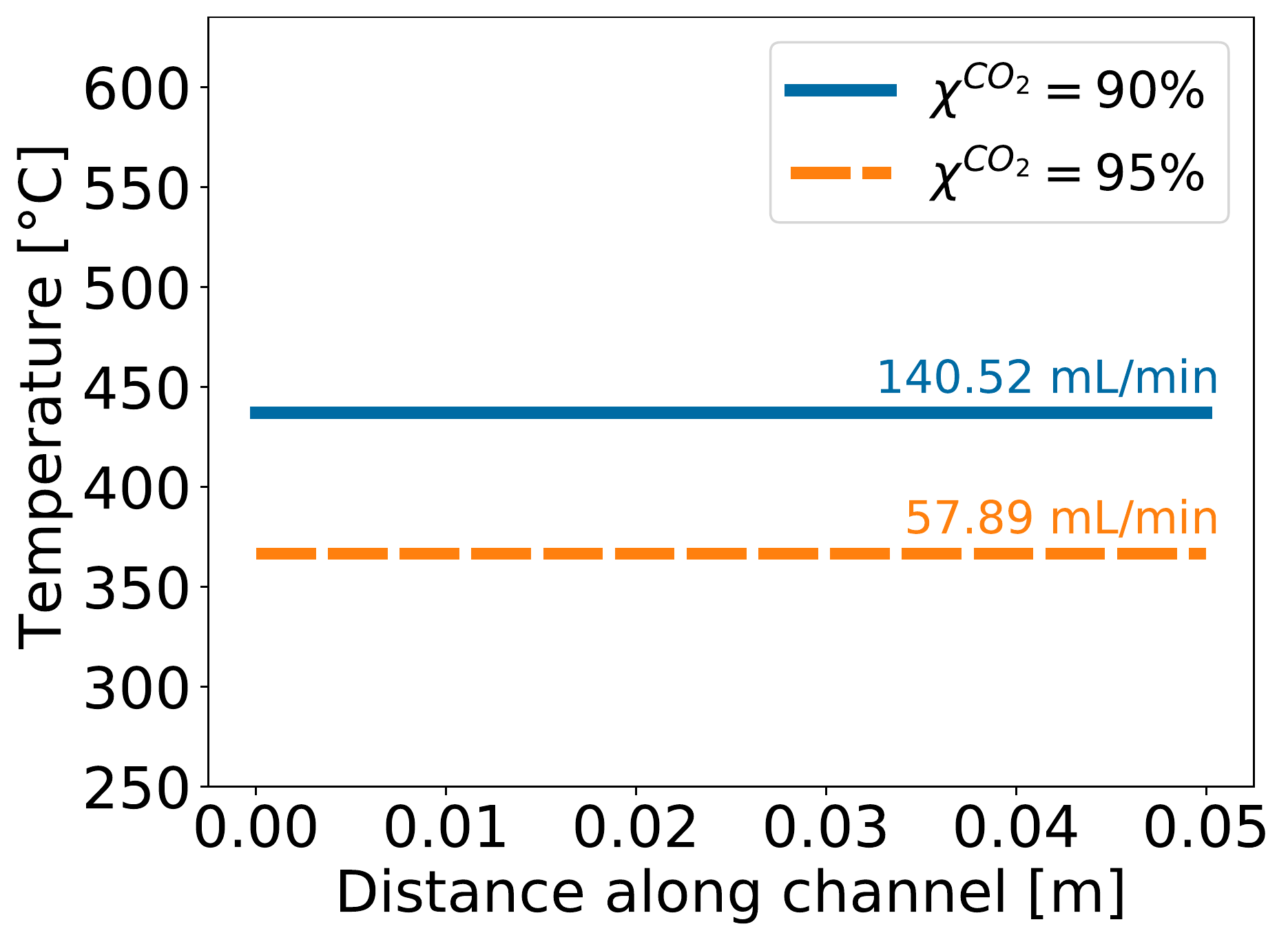}
		\caption{Constant temperature model.}
		%		\caption{\tabular[t]{@{}l@{}}\ce{CO2} conversion: 90\ \%, \\ Flow rate: \SI{140.52}{\milli\liter\per\minute}\endtabular}
	\end{subfigure}
	\hfill
	\begin{subfigure}[t]{0.49\textwidth}
		\centering
		\includegraphics[width=0.9\textwidth]{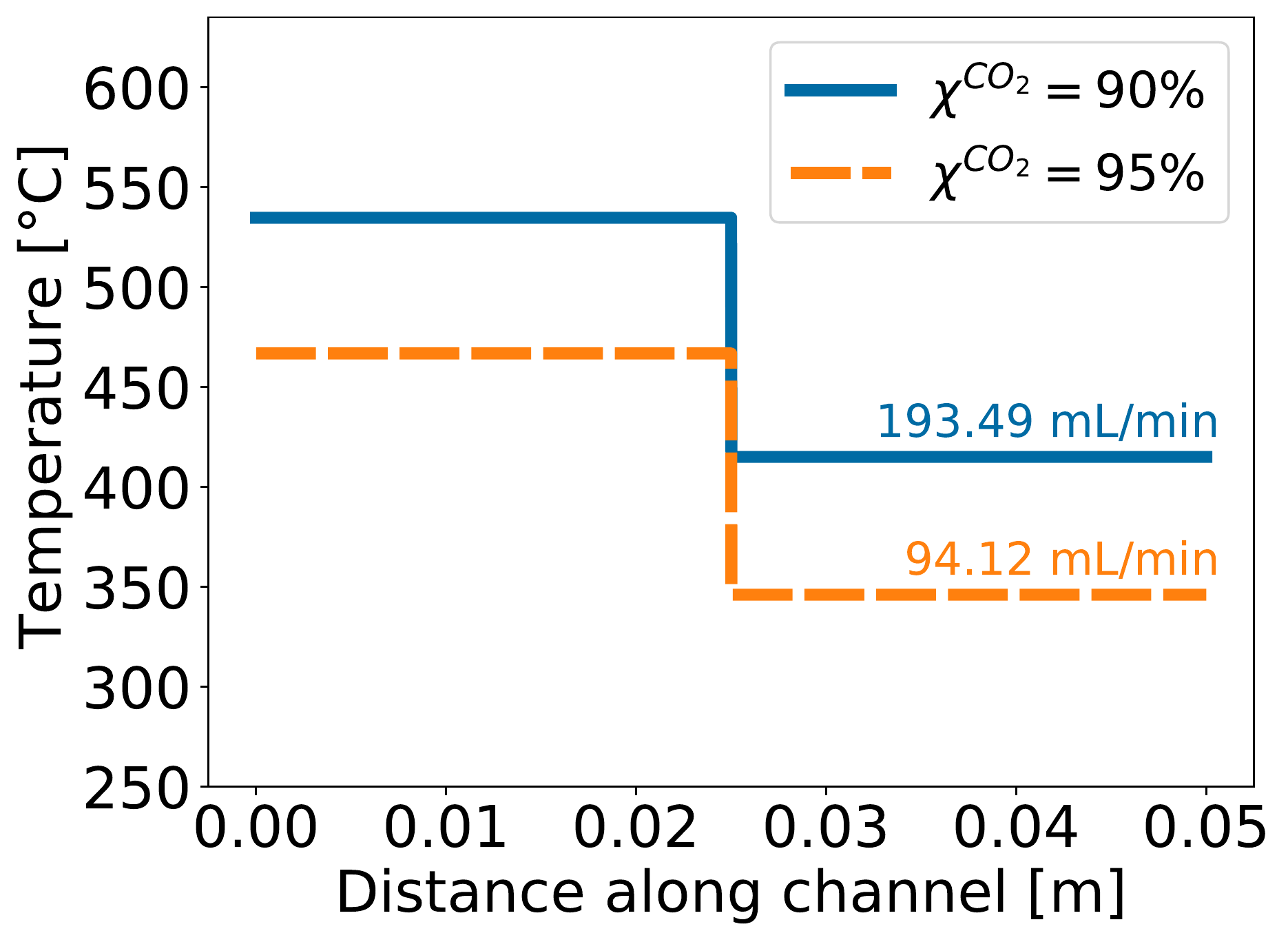}
		\caption{Two temperature stages model.}
		%		\caption{\tabular[t]{@{}l@{}}\ce{CO2} conversion: 90\ \%, \\ Flow rate: \SI{193.49}{\milli\liter\per\minute}\endtabular}
	\end{subfigure}
	\vskip\baselineskip
	\begin{subfigure}[t]{0.49\textwidth}
		\centering
		\includegraphics[width=0.9\textwidth]{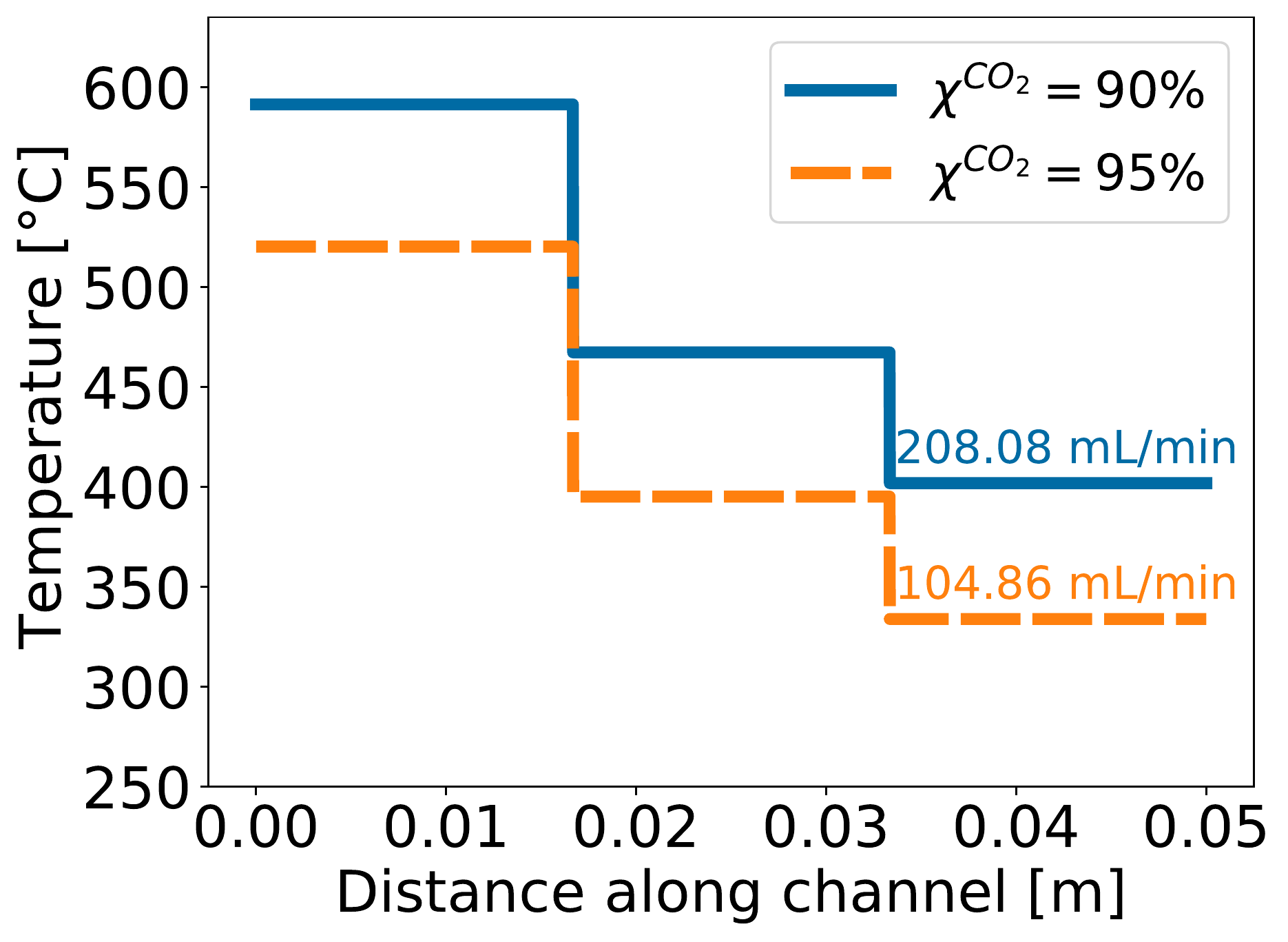}
		\caption{Three temperature stages model.}
		%		\caption{\tabular[t]{@{}l@{}}\ce{CO2} conversion: 90\ \%, \\ Flow rate: \SI{208.08}{\milli\liter\per\minute}\endtabular}
	\end{subfigure}
	\hfill
	\begin{subfigure}[t]{0.49\textwidth}
		\centering
		\includegraphics[width=0.9\textwidth]{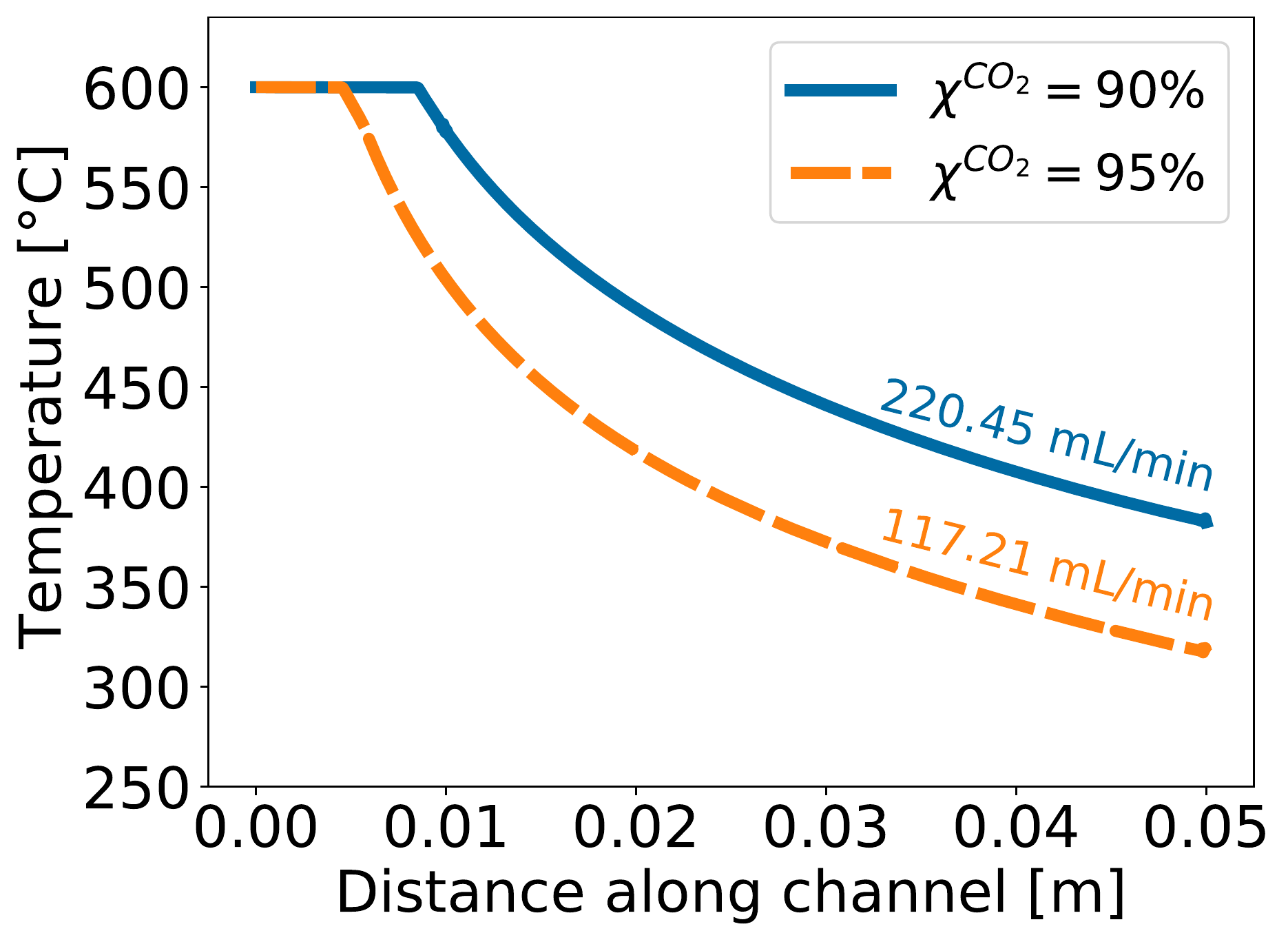}
		\caption{Distributed temperature model.}
		%		\caption{\tabular[t]{@{}l@{}}\ce{CO2} conversion: 90\ \%, \\ Flow rate: \SI{220.45}{\milli\liter\per\minute}\endtabular}
	\end{subfigure}
	\caption{Optimized temperature profiles for problem \eqref{eq:opt_sc_theo} with $\scconst = 0.9$ and $\scconst = 0.95$.}
	\label{fig:opt_sc}
\end{figure}

Let us first investigate the qualitative behavior of the temperature profiles and flow rates, which are depicted in Figure~\ref{fig:opt_sc} for a desired conversion of $\scconst = 0.9$ and $\scconst = 0.95$. We observe that the optimized temperature profiles look very similar to the ones investigated in the previous section (cf. Figure~\ref{fig:opt_cc}). Again, the temperature is higher at the beginning of the reactor and then decreases monotonically along its length for all cases. As expected after our previous investigation, the optimized flow rate increases considerably when using more complex temperature profiles. Moreover, we observe that the largest increase in flow rate happens between the constant and two stage temperature models, whereas increasing the model complexity further has only a smaller influence on the obtained flow rate. For all models the optimized temperature profiles and flow rates are lower for a higher desired conversion $\scconst$. This can be explained through the thermodynamically favorable conditions at lower temperatures and the higher residence times induced by lower flow rates which both lead to an increase in \ce{CO2} conversion (cf. Section~\ref{ssec:sabatier_reaction}).

\begin{table}[!t]
	\centering
	\begin{subtable}{\textwidth}
		\centering
		{\footnotesize
			\rowcolors{2}{\tablegray}{white}
			\setlength{\tabcolsep}{1em}
			\begin{tabular}{S S S S S}
				\toprule
				{\ce{CO2} conversion} & {constant} & {two stages} & {three stages} & {distributed} \\
				\midrule
				\num{85.0}\ \% & 239.74 & 306.03 & 320.89 & 329.42 \\
				\num{87.5}\ \% & 188.01 & 247.08 & 265.08 & 273.80 \\
				\num{90.0}\ \% & 140.52 & 193.49 & 208.08 & 220.45 \\
				\num{92.5}\ \% & 97.20 & 142.71 & 155.83 & 168.65 \\
				\num{95.0}\ \% & 57.89 & 94.12 & 104.86 & 117.21 \\
				\bottomrule
			\end{tabular}
			\caption{Optimized flow rate in [\si{\milli\liter\per\minute}].}
			\label{tab:opt_sc_flow_rate}
		}
	\end{subtable}
	\vskip\baselineskip
	\begin{subtable}{\textwidth}
		\centering
		{\footnotesize
			\rowcolors{2}{\tablegray}{white}
			\setlength{\tabcolsep}{1em}
			\begin{tabular}{c c c c c}
				\toprule
				\ce{CO2} conversion & constant & two stages & three stages & distributed \\
				\midrule
				\num{85.0}\ \% & \num{9.092e-05} & \num{1.161e-04} & \num{1.217e-04} & \num{1.249e-04} \\
				\num{87.5}\ \% & \num{7.340e-05} & \num{9.646e-05} & \num{1.035e-04} & \num{1.069e-04} \\
				\num{90.0}\ \% & \num{5.642e-05} & \num{7.769e-05} & \num{8.355e-05} & \num{8.852e-05} \\
				\num{92.5}\ \% & \num{4.011e-05} & \num{5.889e-05} & \num{6.431e-05} & \num{6.960e-05} \\
				\num{95.0}\ \% & \num{2.454e-05} & \num{3.989e-05} & \num{4.444e-05} & \num{4.968e-05} \\
				\bottomrule
			\end{tabular}
			\caption{Molar product yield in [\si{\mol\per\second}].}
			\label{tab:opt_sc_yield}
		}
	\end{subtable}
	\caption{Numerical results for the optimization problem \eqref{eq:opt_sc_theo}.}
	\label{tab:opt_sc}
\end{table}

The optimized flow rates and absolute product yields for all choices of $\scconst$ are listed in Table~\ref{tab:opt_sc}. Table~\ref{tab:opt_sc_flow_rate} confirms our previous findings, as we see that the optimized flow rate increases with increasing complexity of the temperature models in all cases. As before, we also observe that the biggest difference is between the constant and the two stage model, where the flow rate increases by about \SI{50}{\milli \liter \per \minute}, and additional improvements from the three stage and distributed model are smaller, with about \SI{10}{\milli\liter\per\minute} each. Furthermore, we observe that the optimized flow rate decreases rapidly when the \ce{CO2} conversion is increased. This indicates that higher \ce{CO2} conversions are considerably harder to achieve than lower ones. Considering the total molar product yield, which is depicted in Table~\ref{tab:opt_sc_yield}, we obviously have a similar behavior regarding the complexity of the temperature models as for the flow rate, since the \ce{CO2} conversion restricted by the state constraint for all models. Finally, concerning the influence of the state constraint, it is interesting to observe that the molar product yield decreases with an increase of the desired conversion $\scconst$. In particular, the higher achieved level of \ce{CO2} conversion cannot compensate the necessary reduction in the flow rate, so that overall less product is formed. This implies that one should consider optimizing the flow rate for a desired conversion $\scconst$ as closely as possible to the actual desired product quality since using higher values of $\scconst$ than necessary leads to lower product yields.

\subsection{Realizability of our Results in Practice}
\label{ssec:realizability}

Usually, chemical reactors are equipped with heat exchanger or cooling systems to control the temperature in the reactor. Regarding the realization of the optimized temperature profiles shown in Figures~\ref{fig:opt_cc} and~\ref{fig:opt_sc} with such systems, we have the following remarks. For a microchannel reactor an obvious heat exchanger system is one which is based on microchannels, too, as they have large heat transfer coefficients due to their high specific surface area. Moreover, they already have the appropriate dimensions to be coupled to a microchannel reactor. We refer the reader to, e.g., \cite{khan2011review, naqiuddin2018overview} for an overview of such temperature control systems, as well as to our previous work \cite{Blauth2019Model, Blauth2020Shape}, where we considered the shape optimization of such systems. Naturally, for these systems the question arises, how to align the flow of the coolant with the flow of the chemically reacting gas mixture. Popular choices are given by setups that use co-current, counter-current, or cross-flow configurations \cite{Zanfir2011Optimizing, Park2016Design, Engelbrecht2020Microchannel}. 

Let us start by investigating the counter-current configuration. In this setup, the coolant enters the domain around the area where the gas leaves it. It is at its lowest temperature there since it did not absorb a large amount of heat yet. As the coolant traverses its channel, it is constantly heated up by the reaction and becomes warmest near the reactor inlet. Looking at this from the view of the reactor channel, the gas temperature is at its highest point at the inlet and then gets cooled down progressively as it passes through the reactor. From our results (cf. Figures~\ref{fig:opt_cc} and~\ref{fig:opt_sc}) we can see that this behavior is desirable, so that a suitably chosen counter-current flow configuration could potentially approximate the temperature profiles from the distributed model. These deliberations also directly imply that the co-current flow configuration should be avoided if possible as this would cool the gas right at the inlet and heat it up as it passes through the reactor, achieving the opposite of our goal.

Finally, the cross-flow configuration could also be used to achieve the desired temperatures in the reactor: by using different coolant temperatures in the channels, corresponding to different sections of the reactor, one can achieve high temperatures around the reactor inlet and low temperatures around its outlet. In particular, if a sufficient amount of cooling channels is used, temperature profiles similar to those obtained by the distributed model could potentially be achieved. Moreover, also the piecewise constant temperature profiles could be, approximately, realizable using a cross-flow configuration by dividing the cooling channels into distinct groups, each with the same coolant temperature. In particular, these deliberations show that our choice of the temperature models for the considered optimization problems is sensible as they can be approximately realized in practice by means of appropriate temperature control systems.

%Finally, comparing the obtained molar product yields from the previous section (cf. Table~\ref{tab:opt_cc_yield}), where the flow rate was fixed and we only optimized the wall temperature, with the product yields for problem \eqref{eq:opt_sc_theo} (cf. Table~\ref{tab:opt_sc_yield}), where we optimized the wall temperature and the flow rate, reveals that the optimization of the flow rate can have a significant impact on the overall product yield. Depending on the type of temperature profile, we can increase the product yield by a factor of 4 -- 5 for a desired \ce{CO2} conversion of \num{85}\ \%, and still by up to a factor 2 for \num{95}\ \% conversion. 

In conclusion, we have seen that optimizing a microchannel reactor with methods from PDE constrained optimization shows clear potential for considerable improvements of the reactor's design. In particular, all of our results indicate that optimizing the temperature profile and flow rate of the reactor can increase the obtained product yields considerably. As the considered temperature models are realizable in practice, our results can be readily transferred to the actual application in order to optimize its design. However, we remark that, for this transfer, one also has to consider additional factors, such as energy and cost efficiency, for the particular application at hand.

\section{Conclusion and Outlook}
\label{sec:conclusion}

In this work, we have introduced two mathematical models for chemically reacting fluid flow. First, a three-dimensional model that resolves all quantities over the channel cross section, and second, a one-dimensional model we derived using a homogenization of the reactor over its cross section. We used these to model the Sabatier reaction in a microchannel reactor. The necessary kinetic parameters to compute the reaction rate for these models were determined by means of solving a parameter identification problem constrained by the one-dimensional model using our software package cashocs \cite{Blauth2020CASHOCS} which implements and automates the adjoint approach for solving PDE constrained optimization problems. The obtained results show excellent agreement with the corresponding experimental data published in \cite{Engelbrecht2017Experimentation, Engelbrecht2017Carbon}. We compared both models numerically and found that the one-dimensional model is a very good approximation to the three-dimensional one, with relative errors well below \num{1}\ \%. Finally, we investigated two types of optimal control problems aimed at improving the quality of the reactor. First, we considered the wall temperature of the reactor as optimization variable and maximized the \ce{CO2} conversion for a fixed flow rate. Second, we considered the inlet flow velocity as additional control variable and maximized the mass flow rate of the reactor subject to a state constraint for the \ce{CO2} conversion that ensures a desired quality of the product. The results we obtained for both problems show great potential for optimizing the design of the reactor.

For future work, we propose the extension of our reactor model to include, e.g., more sophisticated reaction mechanisms for the Sabatier process or the RWGS reaction, which would yield a more accurate model. Additionally, it could be worthwhile to investigate other aspects of the reactor, e.g., by considering not only a single channel, but the entire reactor, or by coupling an appropriate temperature control system to the reactor, and to investigate similar optimization problems to the ones in this work in these settings. Another interesting possibility for future research would be to transfer our results into practice and to test them in reality.  Finally, it is also of interest to transfer our setting to other reactions or more complex reactor settings, where, e.g., multiple reactions occur consecutively or in parallel, and to investigate their optimization using similar methods to the ones considered in this paper.

\section*{Acknowledgments}

S. Blauth gratefully acknowledges financial support from the Fraunhofer Institute for Industrial Mathematics ITWM.

%\bibliographystyle{unsrt}
%\bibliography{bibliography.bib}

\appendix
\section{Constitutive Relations}
\label{app:constitutive_relations}

In this appendix we detail the constitutive relations needed to complete the PDE systems \eqref{eq:pde_system} and \eqref{eq:pde1d}.
The first constitutive relation is given by the ideal gas law in its form \eqref{eq:ideal_gas}, i.e.,
\begin{equation*}
	\density = \frac{\pref \avmm}{\gasconst \temperature}.
\end{equation*}
To calculate the average molar mass $\avmm$, we use the relation
\begin{equation*}
	\avmm = \left(\sum_{k=1}^{\nspec} \frac{\massfrac}{\specmm}\right)^{-1},
\end{equation*}
where $\specmm$ is the molar mass of species $k$, which can be found, e.g., in \cite{McBride2002NASA}.

We calculate the specific heat capacity $\hcapacity$, specific enthalpy $\enthalp$, and specific entropy $\entrop$ of the gas mixture as mass averages \cite{Poinsot2005Theoretical} via
\begin{equation*}
	\hcapacity(\massfracvec, \temperature) = \sum_{k=1}^{\nspec} \massfrac \spechcap(\temperature), \qquad \enthalp(\massfracvec, \temperature) = \sum_{k=1}^{\nspec} \massfrac \specenth(\temperature), \qquad \entrop(\massfracvec, \temperature) = \sum_{k=1}^{\nspec} \massfrac \specentr(\temperature),
\end{equation*}
where the corresponding quantities for pure species $\spechcap, \specenth$ and $\specentr$ are determined using the polynomial fits from \cite{McBride2002NASA}
\begin{align*}
	\frac{\spechcap(\temperature)\ \specmm}{\gasconst} &= a_1 \temperature^{-2} + a_2 \temperature^{-1} + a_3 + a_4 \temperature + a_5 \temperature^{2} + a_6 \temperature^{3} + a_7 \temperature^4, \\
	\frac{\specenth(\temperature)\ \specmm}{\gasconst} &= -a_1 \temperature^{-1} + a_2 \log(\temperature) + a_3 \temperature + \frac{a_4}{2} \temperature^2 + \frac{a_5}{3} \temperature^3 + \frac{a_6}{4} \temperature^4 + \frac{a_7}{5} \temperature^5 + b_1, \\
	\frac{\specentr(\temperature)\ \specmm}{\gasconst} &= -\frac{a_1}{2} \temperature^{-2} - a_2 \temperature^{-1} + a_3 \log(\temperature) + a_4 \temperature + \frac{a_5}{2} \temperature^2 + \frac{a_6}{3} \temperature^3 + \frac{a_7}{4} \temperature^4 + b_2.
\end{align*}
This reflects that the specific enthalpy and specific entropy are defined as
\begin{equation*}
	\enthalp(\massfracvec, \temperature) = \int_{\temperature_0}^{\temperature} \hcapacity(\massfracvec, \theta) \dmeas{\theta} + \Delta h_\text{f}^{0}, \quad \text{ and } \quad \entrop(\massfracvec, \temperature) = \int_{\temperature_0}^{\temperature} \frac{\hcapacity(\massfracvec, \theta)}{\theta} \dmeas{\theta} + \Delta s_\text{f}^{0},
\end{equation*}
with specific standard enthalpy of formation $\Delta h_\text{f}^{0}$ and specific standard entropy of formation $\Delta s_\text{f}^{0}$.

For the viscosity $\viscosity$ we have the following. The pure species viscosities $\specvis$ are calculated via the fit given in \cite{Svehla1995Transport}
\begin{equation*}
	\specvis(\temperature) = \exp\left( A_\text{v} \log(\temperature) + B_\text{v} \temperature^{-1} + C_\text{v} \temperature^{-2} + D_\text{v} \right).
\end{equation*}
The mixture-averaged viscosity is then given by \cite{Kee2005Chemically} as
\begin{equation*}
	\viscosity(\massfracvec, \temperature) = \sum_{k=1}^{\nspec} \frac{\molefrac \specvis(\temperature)}{\sum_{j=1}^{\nspec} X_j \Phi_{k, j}},
\end{equation*}
where $\molefrac = \frac{\massfrac\ \avmm}{\specmm}$ is the mole fraction of species $k$, and the weights $\Phi_{k, j}$ are given by
\begin{equation*}
	\Phi_{k, j} = \frac{\left( 1 + \left( \frac{\mu_k(\temperature)}{\mu_j(\temperature)} \right)^{\nicefrac{1}{2}} \left( \frac{M_j}{M_k} \right)^{\nicefrac{1}{4}} \right)^2}{\left( 8 \left( 1 + \frac{\specmm}{M_j} \right) \right)^{\nicefrac{1}{2}}}.
\end{equation*}

Similarly to the viscosity, the thermal conductivity $\conductivity$ of the gas mixture is calculated from the pure species thermal conductivities $\speccond$, which are given by the analogous fit from \cite{Svehla1995Transport}
\begin{equation*}
	\speccond(\temperature) = \exp\left( A_\text{c} \log(\temperature) + B_\text{c} \temperature^{-1} + C_\text{c} \temperature^{-2} + D_\text{c} \right),
\end{equation*}
with the following averaging rule
\begin{equation*}
	\conductivity(\massfracvec, \temperature) = \frac{1}{2} \left( \sum_{k=1}^{\nspec} \molefrac \speccond(\temperature) + \left(\sum_{k=1}^{\nspec} \frac{\molefrac}{\speccond(\temperature)}\right)^{-1}\right).
\end{equation*}

To model the molecular diffusion in the gas, we use the Chapman-Enskog theory and approximate the multicomponent diffusion by mixture-averaged diffusion coefficients which are defined via the following averaging procedure from \cite[Chapter 11]{Kee2005Chemically}
\begin{equation*}
	\diffusioncoeff(\massfracvec, \temperature) = \left( \sum_{\substack{j=1\\ j\neq k}}^{\nspec} \frac{X_j}{D_{k, j}(\temperature)} + \frac{\molefrac}{1 - \massfrac} \sum_{\substack{j=1\\ j\neq k}}^{\nspec} \frac{Y_j}{D_{k, j}(\temperature)} \right)^{-1}.
\end{equation*}
Here, the $D_{k, j}$'s are the binary diffusion coefficients that can be computed via the formulas given, e.g., in \cite{Kee2005Chemically}.

The equilibrium constant $\keqj$ is calculated from thermodynamic relations as follows \cite{Poinsot2005Theoretical, Kee2005Chemically}
\begin{equation*}
	\keqj = \exp\left( \frac{\Delta S_{\text{r}, j}^0(T)}{\gasconst} - \frac{\Delta H_{\text{r}, j}^0(T)}{\gasconst \temperature} \right) \left( \frac{\patm}{\gasconst \temperature} \right)^{\sum_{k=1}^{\nspec} \stoich_{k, j}},
\end{equation*}
where $\patm$ denotes the atmospheric pressure and the argument of the exponential function corresponds to the change in Gibbs free energy for reaction $j$, which is calculated from
\begin{equation*}
	\Delta S_{\text{r}, j}^0(T) = \sum_{k=1}^{\nspec} \stoich_{k, j} \specmm \specentr(T), \quad \text{ and } \quad \Delta H_{\text{r}, j}^0(T) = \sum_{k=1}^{\nspec} \stoich_{k, j} \specmm \specenth(T).
\end{equation*}

%\section{Test}

\end{document}